\documentclass[final,3p,times,twocolumn]{elsarticle}

\usepackage{amsthm}
\usepackage{amsmath}
\usepackage{subfig}
\usepackage{caption}
\usepackage{geometry}
\usepackage{gensymb}
\usepackage{multirow}
\usepackage{algorithmic}
\usepackage[linesnumbered,ruled,vlined]{algorithm2e}

\usepackage{setspace}

\usepackage{gensymb}
\usepackage[utf8]{inputenc}
\usepackage{array}

\usepackage{graphicx}
\usepackage{tikz}
\usepackage{soul}

\usepackage{amssymb}

\usepackage{lineno}

\journal{Computational Geosciences}

\begin{document}

\begin{frontmatter}

\title{Non-intrusive Subdomain POD-TPWL Algorithm for Reservoir History Matching}
\author{Cong Xiao$^{1}$, Olwijn Leeuwenburgh$^{2,3}$, Hai Xiang Lin$^{1}$, Arnold Heemink$^{1}$}

\address{$^{1}$Delft Institute of Applied Mathematics, Delft University of Technology, Mekelweg 4, 2628 CD Delft, the Netherlands}
\address{$^{2}$Civil Engineering and Geosciences, Delft University of Technology, Mekelweg 4, 2628 CD Delft, the Netherlands}
\address{$^{3}$TNO, Princetonlaan 6 PO Box 80015 3508 TA Utrecht, the Netherlands}

\begin{abstract}
This paper presents a non-intrusive subdomain POD-TPWL (SD POD-TPWL) algorithm for reservoir data assimilation through integrating domain decomposition (DD), radial basis function (RBF)
interpolation and the trajectory piecewise linearization (TPWL). It is an efficient approach for model reduction and linearization of general non-linear
time-dependent dynamical systems without intruding the legacy source code. In the  subdomain POD-TPWL algorithm,
firstly, a sequence of snapshots over the entire computational domain are saved and then partitioned into subdomains.
From the local sequence of snapshots over each subdomain, a number of local basis vectors is formed using POD,
and then the RBF interpolation is used to estimate the derivative matrices for each subdomain. Finally,
those derivative matrices are substituted into a POD-TPWL algorithm to form a reduced-order linear model in each subdomain.
This reduced-order linear model makes the implementation of the adjoint easy and resulting in
an efficient adjoint-based parameter estimation procedure. The performance of the new adjoint-based parameter estimation
algorithm has been assessed through several synthetic cases.
Comparisons with the classic finite-difference based history matching show that our proposed subdomain POD-TPWL approach is obtaining comparable results.
The number of full-order model simulations required is roughly 2-3 times the number of uncertain parameters. Using different background parameter realizations,
our approach efficiently generates an ensemble of calibrated models without additional full-order model simulations.
\end{abstract}

\begin{keyword}
Data assimilation \sep reduced-order modeling  \sep model linearization \sep domain decomposition
\end{keyword}

\end{frontmatter}
\section*{Abbreviation}
POD, proper orthogonal decomposition;  RBF, radial basis function;  TPWL, trajectory piecewise linearizaton; DD, domain decomposition; FOM, full-order model

\section{Introduction}
History matching is the process of calibrating uncertain reservoir model parameters such as gridblock permeabilities, porosities, faults multipliers and facies distributions, through minimization of a cost function that quantifies the misfit between simulated and observed data (typically well data such as oil or water rates or bottom-hole pressure, but possibly also 4D seismic data. If the gradient of the cost function with respect to parameters can be computed using the adjoint of the reservoir model, history matching problems can be efficiently solved using a gradient-based minimization algorithm \cite{courant1962methods}.  In general, significant effort is required to obtain and maintain a correct implementation of the adjoint model for complex nonlinear simulation models. Such implementations are generally instrusive, that is, they require access to the model code, which may not always be possible.

Many efforts have been taken to make the implementation of the adjoint model more feasible. One way is to replace the original complex model with a surrogate so that the construction of the adjoint model becomes easier. Courtier et al (1994) \cite{courtier1994strategy} proposed an incremental approach by replacing
a high resolution nonlinear model with an approximated linear model so that the adjoint model can be more easily obtained. Liu et al (2008, 2009) \cite{liu2008ensemble}, \cite{liu2009ensemble} developed an ensemble-based four-dimensional variational (En4DVar) data assimilation scheme where the approximated linear model is constructed using an ensemble of model forecasts. Recently, to extend the ensemble-based tangent linear model (TLM) to more realistic applications, Frolov and Bishop et al (2016, 2017) \cite{frolov2016localized}, \cite{bishop2017local} incorporated a local
ensemble tangent linear model (LETLM) into 4D-Var scheme. The LETLM has the ability to capture localized physical features of dynamic models with relatively small ensemble size. However, the construction of a tangent linear model becomes intractable for high-dimensional systems. Proper Orthogonal Decomposition (POD), a model order reduction method, is a possible approach to decrease the dimensionality of the original model. The POD approach has been applied to various disciplines, including reservoir model simulations \cite{heijn2003generation}, \cite{markovinovic2006accelerating} and has in some cases shown significant speed up \cite{cardoso2009development} .

The combination of model linearization and model reduction techniques has the potential to further ease the implementation of adjoint models for high-dimensional complex dynamic systems. Vermeulen and Heemink (2006) \cite{vermeulen2006model} combined POD and a non-intrusive perturbation-based linearization method to build a reduced-order linear approximation of the original high-dimensional non-linear model. The adjoint of this reduced-order linear model can be easily constructed and therefore the minimization of the objective function can be handled efficiently. Altaf et al (2009) \cite{altaf2009inverse} and Kaleta et al (2011) \cite{kaleta2011model} applied this method to a coastal engineering and reservoir history matching problem, respectively.

Alternatively, Trajectory Piecewise Linearization (TPWL) can be classified
as a model-intrusive linearization method. In TPWL, a number of full-order 'training' runs is first simulated, and then a linear model
is generated through first-order expansion around the 'closest' training trajectories. In reservoir engineering, Cardoso et al (2010) \cite{cardoso2010linearized} was the first to integrate the POD and the TPWL methods and applied this strategy to oil production optimization. He et al. (2013, 2014) applied the POD-TPWL method to both reservoir history matching and production optimization \cite{he2013reduced}, \cite{he2014reduced}. These studies suggested that POD-TPWL has the potential to significantly reduce the runtime for subsurface flow problems \cite{he2011enhanced}. A drawback, however, is that the POD-TPWL method requires access to derivative matrices used internally by the numerical solver, and
therefore cannot be used with most commercial simulators \cite{he2013reduced},\cite{trehan2016trajectory}. And although the traditional construction of a reduced-order linear model \cite{vermeulen2006model}, \cite{altaf2009inverse}, \cite{kaleta2011model} is model non-intrusive,
the required derivative information is estimated using a global perturbation-based finite difference method, which needs a large number of full-order simulations and is therefore computationally demanding. Furthermore, the global perturbation also hinders the extension of this method to large-scale reservoir history matching which requires retaining many POD patterns. In order to avoid model intrusion and numerous full-order simulations, we propose to incorporate domain decomposition (DD) and radial basis function (RBF) interpolation into POD-TPWL to develop a new non-intrusive subdomain POD-TPWL algorithm.

RBF interpolation is mainly used to construct surrogate models, and has been applied e.g. to reservoir engineering and fluid dynamics \cite{klie2013unlocking}, \cite{xiao2016nona}, \cite{xiao2016nonb}. Recently,  Bruyelle et al (2014) \cite{bruyelle2014neural} applied the neural network-based RBF to obtain the first-order and second-order derivative information of a reservoir model and estimate the gradients and Hessian matrix for reservoir production optimization. The accuracy of RBF-based gradient approximation is determined by the sampling strategies of the interpolation data \cite{bruyelle2014neural}. For high dimensional problems, the classical global RBF interpolation algorithm requires a large number of interpolation data to capture the flow dynamic as much as possible \cite{chinchapatnam2007domain}. Moreover, the global RBF algorithm can cause some spurious long-distance correlations, which implies the possibilities to avoid some redundant interpolation data. This motivates us to develop a subdomian RBF interpolation technique for reservoir models where the domain decomposition (DD) technique potentially allows us to apply the methodology to large-scale problems. Different local RBF interpolation schemes are considered based on the details of local flow dynamics in each subdomain.
The domain decomposition technique first introduced in the work of Przemieniecki \cite{przemieniecki1963matrix} has been applied to various fields \cite{bian2015multi}. Lucia et al. (2003) \cite{lucia2003reduced} first introduced the DD method to model-order reduction for accurately tracking a moving strong shock wave. Subsequently, the DD method has also been applied to non-linear model reduction problems \cite{baiges2013domain},  \cite{amsallem2012nonlinear}, \cite{chaturantabut2017temporal}.

This paper presents a new non-intrusive subdomain POD-TPWL algorithm for subsurface flow problem.
The key idea behind this subdomain POD-TPWL is to integrate the DD method and RBF algorithm into a model linearization technique based on the POD-TPWL.
After constructing the reduced-order linear model using subdomain POD-TPWL algorithm, because of the linearity in the reduced-order subspace,
the implementation of adjoint model is easy and, thus, it is convenient to incorporate this reduced-order linear model into a gradient-based reservoir history matching procedure.
The runtime speedup and the robustness of the new history matching algorithm have been assessed through several synthetic cases.

This paper is arranged as follows: The history matching problem and the classical adjoint-based solution approach are described in Section 2.
Section 3 contains the mathematical background of the traditional POD-TPWL. Section 4 gives the mathematical descriptions of domain decomposition (DD) and radial basis function (RBF) interpolation,
which are used to develop the non-intrusive subdomain POD-TPWL algorithm. In addition, a workflow for combining subdomain POD-TPWL with an adjoint-based history matching algorithm is described.
Section 5 discusses and evaluates an application of the new history matching workflow to some numerical 'twin' experiments involving synthetic reservoir models.
Finally, Section 6 summarizes our contribution and discusses future work.

\section{Problem Description}
A single simulation step of a discretized two-phase oil-water reservoir system is described as follows,
\begin{equation}
\label{eq1}
\textbf{x}^{n+1} = \textbf{f}^{n+1}(\textbf{x}^{n},\boldsymbol{\beta}), \quad n=1,\cdot\cdot\cdot,N
\end{equation}
where the dynamic operator ${\textbf{f}}^{n+1}$: $R^{2N_{g}}{\rightarrow}R^{2N_{g}}$ represents the nonlinear time-dependent model evolution, $ \textbf{x}^{n+1}{\in}R^{2N_{g}} $ represents the state vector (pressure and saturation in every gridblock), $N_{g}$ is the total number of gridblocks,
${\textit{n}}$ and $ {\textit{n}}+1 $ indicate the timesteps, $ \textit{N} $ denotes the total number of simulation steps, and $\boldsymbol{\beta}$ denotes the vector of uncertain parameters, which is the spatial permeability field in our case. For more details about the discretization of the governing equations, see e.g., \cite{peacemanfundamentals}.

The relationship between simulated data $ \textbf{y}^{m+1} $ and state vector $ \textbf{x}^{m+1} $ can be described by a nonlinear operator ${\textbf{h}}^{m+1}$: $R^{2N_{g}}{\rightarrow}R^{N_{d}}$, which, in our case, represents the well model (for seismic data another model would be needed). $ N_{d}$ is the number of measurements at each timestep. The simulated measurements are therefore described by
\begin{equation}
\label{eq2}
\textbf{y}^{m+1}=\textbf{h}^{m+1}(\textbf{x}^{m+1},\boldsymbol{\beta}), \quad m=1,\cdot\cdot\cdot,N_{0}
\end{equation}
where $\textit{N}_{0}$ is the number of timesteps at which measurements are taken.

The history matching process calibrates the uncertain parameters by minimizing a cost function defined as a sum of weighted squared differences between observed and modeled measurements (data). Additional incorporation of prior information into the cost function as a regularization term can further constrain the minimization procedure and make the history matching problem well-posed \cite{oliver2008inverse}. Eventually, the cost function is described by the sum of two terms.
\begin{align}
\label{eq3}
J(\textbf{x}^{1},\cdot\cdot\cdot,\textbf{x}^{n},\cdot\cdot\cdot,\textbf{x}^{N},\boldsymbol{\beta}) = \frac{1}{2}(\boldsymbol{\beta}-\boldsymbol{\beta}_{p})^{T}{\textbf{R}_{p}}^{-1}(\boldsymbol{\beta}-\boldsymbol{\beta}_{p}) \notag \\
+\frac{1}{2}\sum_{m=1}^{N_{0}}(\textbf{d}_{o}^{m}-\textbf{h}^{m}(\textbf{x}^{m},\boldsymbol{\beta}))^{T}{\textbf{R}_{m}}^{-1}(\textbf{d}_{o}^{m}-\textbf{h}^{m}(\textbf{x}^{m},\boldsymbol{\beta}))
\end{align}
where $\textbf{d}_{o}^{m}$ represents the vector of observed data at timestep $\textit{m}$.

In twin experiments $\textbf{d}_{o}^{m}$ is generated by adding some noise, e.g $\textbf{r}^{m}$, to the data $\textbf{y}_t^{m}$ simulated with a 'truth' model. We will assume here that $\textbf{r}^{m}$ is a time-dependent vector of observation errors at time level ${\textit{m}}$, which is uncorrelated over time, and satisfies the Gaussian distribution $G (\textbf{0},\textbf{R}_{m})$ where $\textbf{R}_{m}$ is the observation error covariance matrix at the timestep $\textit{m}$. $\boldsymbol{\beta}_{p}$ represents the prior parameter vector, and $\textbf{R}_{p}$ represents the error covariance matrix of the prior parameters, which characterizes the uncertainty in the prior model. A gradient-based optimization algorithm can be used to determine a parameter set that is not too far away from the prior information, while minimizing the misfit between the observed and simulated data.

The key step of a gradient-based minimization algorithm is to determine the gradient of the cost function with respect to the parameters. The gradient of the cost function can be formulated by introducing the adjoint model as follows (more details about the mathematical derivation can be found in \cite{jansen2011adjoint}),

\begin{align}
\label{eq4}
[\frac{d J}{d \boldsymbol{\beta}}]^{T} & = {\textbf{R}_{p}}^{-1}(\boldsymbol{\beta}-\boldsymbol{\beta}_{p})-\sum_{n=1}^{N}[\boldsymbol{\lambda}^{n}]^{T}{\frac{\partial \textbf{f}^{n}}{\partial \boldsymbol{\beta}}} \notag \\
& -\sum_{m=1}^{N_{0}}[\frac{\partial \textbf{h}^{m}(\textbf{x}^{m},\boldsymbol{\beta})}{\partial \boldsymbol{\beta}}]^{T}{\textbf{R}_{m}}^{-1}(\textbf{d}_{o}^{m}-\textbf{h}^{m}(\textbf{x}^{m},\boldsymbol{\beta}))
\end{align}
where the adjoint model in terms of the Lagrange multipliers $\boldsymbol{\lambda}^{n}$ is given by
\begin{equation}
\label{eq5}
\boldsymbol{\lambda}^{n}=[\frac{\partial \textbf{f}^{n+1}}{\partial \textbf{x}^{n}}]\boldsymbol{\lambda}^{n+1}+[\frac{\partial \textbf{h}^{n}(\textbf{x}^{n},\boldsymbol{\beta})}{\partial \textbf{x}^{n}}]^{T}{\textbf{R}_{n}}^{-1}(\textbf{d}_{o}^{n}-\textbf{h}^{n}(\textbf{x}^{n},\boldsymbol{\beta}))
\end{equation}
for $\textit{n}$ = $\textit{N},\cdot\cdot\cdot\,1 $ with an ending condition $\boldsymbol{\lambda}^{N+1}=0$. This adjoint approach has a high computational efficiency because just one forward simulation and one backward simulation are required to compute the gradient, independent on the size of the variable vector.
It should be pointed out that four derivative terms, e.g, $\frac{\partial \textbf{h}^{m}(\textbf{x}^{m},\boldsymbol{\beta})}{\partial \boldsymbol{\beta}}$, $\frac{\partial \textbf{h}^{n}(\textbf{x}^{n},\boldsymbol{\beta})}{\partial \textbf{x}^{n}}$, $\frac{\partial \textbf{f}^{n}}{\partial \boldsymbol{\beta}}$ and $\frac{\partial \textbf{f}^{n+1}}{\partial \textbf{x}^{n}}$, are required in the adjoint method. We will give detailed descriptions of how to efficiently obtain these four terms using our proposed subdomain POD-TPWL algorithm in the following sections.

\section{POD-TPWL algorithm}
In the TPWL scheme, one or more full order "training" runs using a
set of perturbed parameters are simulated first. The states and the derivative information at each time step from these runs are used to construct the TPWL surrogate. Given the state $\textbf{x}^{n} $ and parameters $\boldsymbol{\beta}$, the state $ \textbf{x}^{n+1} $ is approximated as a first-order expansion around the training solution
$(\textbf{x}_{tr}^{n+1} ,\textbf{x}_{tr}^{n},\boldsymbol{\beta}_{tr})$ as follows,
\begin{equation}
\label{eq6}
\textbf{x}^{n+1} {\approx} \textbf{x}_{tr}^{n+1}+\textbf{E}^{n+1}(\textbf{x}^{n}-\textbf{x}_{tr}^{n})+\textbf{G}^{n+1}(\boldsymbol{\beta}-\boldsymbol{\beta}_{tr})
\end{equation}
\begin{equation}
\label{eq7}
\textbf{E}^{n+1}=\frac{\partial \textbf{f}^{n+1}}{\partial \textbf{x}_{tr}^{n}}, \quad \textbf{G}^{n+1}=\frac{\partial \textbf{f}^{n+1}}{\partial \boldsymbol{\beta}_{tr}}
\end{equation}

The training solution $(\textbf{x}_{tr}^{n+1},\textbf{x}_{tr}^{n},\boldsymbol{\beta}_{tr})$ is chosen to be as 'close' as possible to the state $\textbf{x}^{n}$. A detailed description of the criterion for closeness can be found in \cite{he2015constraint}. The matrices $\textbf{E}^{n+1}\in R^{2N_{g} \times 2N_{g}}$ and $\textbf{G}^{n+1}\in R^{2N_{g} \times N_{g}}$ represent the derivative of the dynamic model (Eq.\ref{eq1}) at timestep $\textit{n}$+1 with respect to states $\textbf{x}_{tr}^{n}$ and parameters $\boldsymbol{\beta}_{tr}$ respectively. Eq.\ref{eq6} is, however, still in a high-dimensional space, e.g, $\textbf{x}^{n+1} \in R^{2N_{g}}$, and $\boldsymbol{\beta}\in R^{N_{g}}$, which motivates the development of the POD-TPWL algorithm \cite{he2015constraint}.

Proper Orthogonal Decomposition (POD) provides a means to project the high-dimensional states into an optimal lower-dimensional subspace. The basis of this subspace is obtained by performing a Singular Value Decomposition (SVD) of a snapshot matrix containing the solution states at selected time steps (snapshots) computed from training simulations. The state vector $\textbf{x}$ can then be  represented in terms of the product of a coefficient
vector $\boldsymbol{\psi}$ and a matrix of basis vectors $\boldsymbol{\phi}$
\begin{equation}
\label{eq8}
\textbf{x} = \boldsymbol{\phi} \boldsymbol{\psi}
\end{equation}
Let $\boldsymbol{\phi}_{p}$ and $\boldsymbol{\phi}_{s}$ represent separate matrices of basis vectors for pressure and saturation respectively. In general there is no need to contain all columns of the left singular matrix in $\boldsymbol{\phi}_{p}$ and $\boldsymbol{\phi}_{s}$ (see e.g. \cite{he2015constraint}) and a reduced state vector representation can be obtained by selecting only the first columns according to e.g. an energy criterion. To normalize the reduced state vector, the columns of $\boldsymbol{\phi}_{p}$ are
determined by multiplying the left singular matrix $ \textbf{U}_{p}$ with the singular value matrix $\boldsymbol{\Sigma}_{p}$ (and similarly for saturation), i.e.

\begin{equation}
\label{eq9}
\boldsymbol{\phi}_{p} = \textbf{U}_{p} \boldsymbol{\Sigma}_{p}, \quad  \boldsymbol{\phi}_{s} = \textbf{U}_{s} \boldsymbol{\Sigma}_{s} \quad .
\end{equation}

In this paper, we use Karhunen-Loeve expansion (KLE) to parameterize the parameter space. KLE
reduces the dimension of the parameter vector by projecting the high-dimensional parameter into an optimal lower-dimensional subspace \cite{fukunaga1970application}. The basis of this subspace is obtained by performing an eigenvalue decomposition of the prior parameter covariance matrix $\textbf{R}_{p}$.
If this covariance matrix is not accessible the basis can alternatively be obtained from an SVD decomposition of matrix holding an ensemble of prior parameter realizations with ensemble mean $\boldsymbol{\beta}_{b} = \boldsymbol{\beta}_{p}$. Including normalization of the reduced parameter vector, a random parameter vector sample $\boldsymbol{\beta}$ can be generated as follows,
\begin{equation}
\label{eq10}
\boldsymbol{\beta} = \boldsymbol{\beta}_{b}+ \boldsymbol{\phi}_{\boldsymbol{\beta}} \boldsymbol{\xi}, \quad \text{with}\quad \boldsymbol{\phi}_{\boldsymbol{\beta}} = \textbf{U}_{\boldsymbol{\beta}} \boldsymbol{\Sigma}_{\beta}
\end{equation}
where $\boldsymbol{\phi}_{\boldsymbol{\beta}}$ denotes the matrix of parameter basis vectors,  $\textbf{U}_{\boldsymbol{\beta}}$ and $\boldsymbol{\Sigma}_{\beta}$ are the left singular matrix and singular value matrix of the parameter matrix respectively, and $\boldsymbol{\xi}$ denotes a vector with independent Gaussian random variables with zeros mean and unit variance. A reduced parameter space representation can again be obtained by selecting only the first several columns of $\boldsymbol{\phi}_{\boldsymbol{\beta}}$ according to e.g. an energy criterion.

The number of the retained columns for basis matrix (denoted as $l_{p}$ and $l_{s}$ for pressure and saturation, $l_{\boldsymbol{\beta}}$ for parameter, respectively) is determined through an energy criterion \cite{he2015constraint}. We take $\boldsymbol{\phi}_{p}$ as an example. We first compute the total energy $E_{t}$, which is defined as $E_{t}=\sum_{i=1}^{L}\nu_{i}^{2}$, where $\nu_{i}$ denotes the $\textit{i}$-th singular value of snapshot matrix for pressure. The energy associated with the first $l_{p}$ singular vectors is given by $E_{l_{p}}=\sum_{i=1}^{l_{p}}{\nu_{i}}^{2}$. Then $l_{p}$ is determined such that $E_{l_{p}}$exceeds a specific fraction of $E_{t}$. The same process can be assigned to determine $l_{s}$ and $l_{\boldsymbol{\beta}}$.

Substituting Eq.\ref{eq8} and Eq.\ref{eq10} into Eq.\ref{eq6}, we obtain the following POD-TPWL formula,
\begin{equation}
\label{eq11}
\boldsymbol{\psi}^{n+1} {\approx} \boldsymbol{\psi}_{tr}^{n+1}+\textbf{E}_{\boldsymbol{\psi}}^{n+1}(\boldsymbol{\psi}^{n}-\boldsymbol{\psi}_{tr}^{n})+\textbf{G}_{\boldsymbol{\xi}}^{n+1}(\boldsymbol{\xi}-\boldsymbol{\xi}_{tr})
\end{equation}
\begin{equation}
\label{eq12}
\textbf{E}_{\boldsymbol{\psi}}^{n+1}=\boldsymbol{\phi}^{T}\frac{\partial \textbf{f}^{n+1}}{\partial \textbf{x}_{tr}^{n}}\boldsymbol{\phi}, \quad \textbf{G}_{\boldsymbol{\xi}}^{n+1}=\boldsymbol{\phi}^{T}\frac{\partial \textbf{f}^{n+1}}{\partial \boldsymbol{\beta}_{tr}}\boldsymbol{\phi}_\beta
\end{equation}

Similarly, the well model as Eq.\ref{eq2} is also linearized around a close training solution $(\boldsymbol{\psi}_{tr}^{n+1}, \boldsymbol{\xi}_{tr})$ in the reduced space as follows,
\begin{equation}
\label{eq13}
\textbf{y}^{m+1} {\approx} \textbf{y}_{tr}^{m+1}+\textbf{A}_{\boldsymbol{\psi}}^{m+1}(\boldsymbol{\psi}^{m+1}-\boldsymbol{\psi}_{tr}^{m+1})+\textbf{B}_{\boldsymbol{\xi}}^{m+1}(\boldsymbol{\xi}-\boldsymbol{\xi}_{tr})
\end{equation}
\begin{equation}
\label{eq14}
\textbf{A}_{\boldsymbol{\psi}}^{m+1}=\frac{\partial \textbf{h}^{m+1}}{\partial \textbf{x}_{tr}^{m+1}}\boldsymbol{\phi}, \quad \textbf{B}_{\boldsymbol{\xi}}^{m+1}=\frac{\partial \textbf{h}^{m+1}}{\partial \boldsymbol{\beta}_{tr}}\boldsymbol{\phi}_{\boldsymbol{\beta}}
\end{equation}

Eq.\ref{eq11} and Eq.\ref{eq13} represent the POD-TPWL system for reservoir model and well model in the reduced-order space, respectively. In general the traditional POD-TPWL method modifies the source code to output all derivative matrices \cite{he2015constraint}. In this paper, we integrate domain decomposition technique and radial basis function interpolation to
approximately estimate these derivative matrices without accessing to the code. These derivative matrices then are substituted into POD-TPWL algorithm to form a subdomain reduced-order linear model.

\section{Adjoint-based history matching using subdomain reduced-order linear model}
This section describes the mathematical background of domain decomposition (DD), and radial basis function (RBF) interpolation, which are used to construct the subdomain non-intrusive reduced
order linear model. In addition, how to incorporate this reduced-order linear model into the adjoint-based history matching is described in the last subsection.

\subsection{Domain Decomposition Method}
A 2D or 3D computational domain is denoted as $\Omega$. The entire  domain $\Omega$ is assumed to be decomposed into $\textit{S}$ non-overlapping subdomains $\Omega^{d}$ , $d \in \{1,2,\cdot\cdot\cdot,S\}$
(such that $\Omega=\bigcup_{d=1}^{S}\Omega^{d} $ and $ \Omega^{i} \cap \Omega^{j}=0 $ for $ i\neq j$) and each subdomain has local unknowns.e.g, local pressure and saturation variables.
In each subdomain $\Omega^{d}$, the generated global snapshots within that subdomain are used to construct a set of local POD basis functions $\boldsymbol{\phi}^{d}$ and the corresponding POD coefficients $\boldsymbol{\psi}^{d,n+1}$ at the timestep $\textit{n}$+1 as described
in the previous section. For each subdomain $\Omega^{d}$, the reservoir dynamic model as Eq.\ref{eq1} is modified to represent the underlying dynamic system associated with this subdomain $\Omega^{d}$ and its surrounding subdomains $\Omega^{sd}$ in the reduced subspace, and can be reformulated as
\begin{equation}
\label{eq15}
\boldsymbol{\psi}^{d,n+1} = \boldsymbol{\pounds}^{d,n+1}(\boldsymbol{\psi}^{d,n},\boldsymbol{\psi}^{sd,n+1},\boldsymbol{\xi})
\end{equation}

\begin{figure}[!h]
\centering\includegraphics[width=0.6\linewidth]{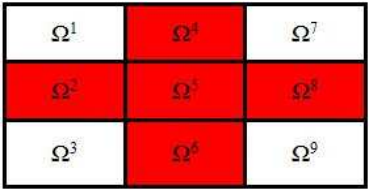}
\caption{Illustration of domain decomposition in a 2-D case}\label{fig1}
\end{figure}

The well model represents the underlying dynamic system just associated with this subdomain $\Omega^{d}$, and can be given by
\begin{equation}
\label{eq16}
\textbf{y}^{d,m+1} = \boldsymbol{\hbar}^{d,n+1}(\boldsymbol{\psi}^{d,m+1},\boldsymbol{\xi})
\end{equation}

where, vector $\boldsymbol{\psi}^{d,n}$ denotes the set of POD coefficients at the time level $\textit{n}$ for the subdomain $\Omega^{d}$, $\boldsymbol{\psi}^{sd,n+1}$ denotes the set of POD coefficients at time level $\textit{n}$+1 for the surrounding subdomains $\Omega^{sd}$. In a 2-D case, the number of surrounding subdomains
associated to this subdomain $\Omega^{d}$ is between 2 and 4, see a simple example in Fig.\ref{fig1}, which shows a maximum of four surrounding subdomains connected with the subdomain $\Omega^{5}$, three surrounding subdomains connected with the subdomain $\Omega^{2}, \Omega^{4}, \Omega^{6}, \Omega^{8}$, and two surrounding subdomains connected with the subdomain $\Omega^{1}, \Omega^{3},
\Omega^{7}, \Omega^{9}$.

The POD-TPWL algorithm forms a reduced-order tangent linear
model (TLM) of the original nonlinear model. Recently, a local ensemble tangent linear model (LETLM) has been developed to capture localized physical features of dynamic models with relatively small ensemble size \cite{frolov2016localized}, \cite{bishop2017local}. The key point of LETLM is to identify the influence area which is very
similar to the purpose of domain decomposition described here. However, this LETLM needs to sequentially construct the TLM piecewise for each state variable, which results in a large number of full-order model simulations and
overwhelming programming efforts. In this study, we construct the reduced-order tangent linear model (TLM) piecewise for each subdomain instead of each state variable. We propose to use RBF interpolation to obtain the
derivative matrices that are required by the POD-TPWL. In addition, domain decomposition has the abilities to
efficiently capture localized physical features \cite{chinchapatnam2007domain}, and therefore has the potential to improve the derivative estimate by local low-dimensional RBF interpolation which will be described in the next subsections.

\subsection{Radial Basis Function Interpolation}
RBF interpolation can be classified as a data-driven interpolation method, which is mainly used to construct surrogate model \cite{klie2013unlocking}, \cite{xiao2016nona}. High-dimensional interpolation needs a large number of data to obtain a satisfactory accuracy, a phenomenon often referred to as the
\textquotedblleft curse of dimensionality\textquotedblright. To remedy this difficulty, domain decomposition technique approximates the global domain by the sum of the local subdomains, and therefore can be applied to form a locally low-dimensional RBF interpolation.

For subdomain $\Omega^{d}$, let $\boldsymbol{\pounds}^{d,n+1}(\boldsymbol{\psi}^{d,n},\boldsymbol{\psi}^{sd,n+1},\boldsymbol{\xi})$ denote a RBF interpolation function for the POD coefficient $\boldsymbol{\psi}^{d,n+1}$
at the time level $\textit{n}$+1. The RBF interpolation function is
a linear combination of $\textit{M}$ radial basis functions in the form of,
\begin{align}
\label{eq17}
& \boldsymbol{\pounds}^{d,n+1}(\boldsymbol{\psi}^{d,n},\boldsymbol{\psi}^{sd,n+1},\boldsymbol{\xi})
= \notag \\
& \sum_{j=1}^{M}\boldsymbol{\omega}_{j}^{d,n+1}\times \theta (||(\boldsymbol{\psi}^{d,n},\boldsymbol{\psi}^{sd,n+1},\boldsymbol{\xi})-(\boldsymbol{\psi}_{j}^{d,n},\boldsymbol{\psi}_{j}^{sd,n+1},\boldsymbol{\xi}_{j})||)
\end{align}
where, $\boldsymbol{\omega}^{d,n+1} $ is a weighting coefficient vector of size $\textit{M}$ (number of training runs). $||(\boldsymbol{\psi}^{d,n},\boldsymbol{\psi}^{sd,n+1},\boldsymbol{\xi})-(\boldsymbol{\psi}_{j}^{d,n},\boldsymbol{\psi}_{j}^{sd,n+1},\boldsymbol{\xi}_{j})||$ is a scalar distance using $\textit{L}_{2} $ norm. $\theta$ is a set of specific radial basis functions.

The specific coefficient $\boldsymbol{\omega}_{j}^{d,n+1}$ is determined so as to ensure that the interpolation function values $\boldsymbol{\pounds}^{d,n+1}$ at the training data points
$(\boldsymbol{\psi}_{j}^{d,n},\boldsymbol{\psi}_{j}^{sd,n+1},\boldsymbol{\xi}_{j})$, matches the given data  $\boldsymbol{\psi}_{j}^{d,n+1}$ exactly. This can be expressed by,
\begin{equation}
\label{eq18}
\textbf{D}^{d,n+1}\boldsymbol{\omega}^{d,n+1}=\textbf{Z}^{d,n+1}
\end{equation}
where
\begin{align}
\label{eq19}
\centering
& \textbf{D}^{d,n+1}=\begin{bmatrix} \theta(l^{n+1}(1,1)) & ... & \theta(l^{n+1}(1,M)) \\ . & \theta(l^{n+1}(i,j) & . \\ \theta(l^{n+1}(M,1)) & ... & \theta(l^{n+1}(M,M)) \end{bmatrix} \notag \\
& l^{n+1}(i,j)=||(\boldsymbol{\psi}_{i}^{d,n},\boldsymbol{\psi}_{i}^{sd,n+1},\boldsymbol{\xi}_{i})-(\boldsymbol{\psi}_{j}^{d,n},\boldsymbol{\psi}_{j}^{sd,n+1},\boldsymbol{\xi}_{j})||, \notag \\
& i=1,\cdot\cdot\cdot,M ; j=1,\cdot\cdot\cdot,M
\end{align}
\begin{equation}
\label{eq20}
\boldsymbol{\omega}^{d,n+1}=[\boldsymbol{\omega}_{1}^{d,n+1}, \boldsymbol{\omega}_{2}^{d,n+1}, \cdot\cdot\cdot, \boldsymbol{\omega}_{M}^{d,n+1}]^{T}
\end{equation}
\begin{equation}
\label{eq21}
\textbf{Z}^{d,n+1}=[\boldsymbol{\psi}_{1}^{d,n+1}, \boldsymbol{\psi}_{2}^{d,n+1}, \cdot\cdot\cdot, \boldsymbol{\psi}_{M}^{d,n+1}]^{T}
\end{equation}

The weighting coefficients are determined by solving the linear system of equation Eq.\ref{eq18}. A list of well-known radial basis function are provided in Table \ref{tab1}.
In general, some different type of radial basis function $\theta$ can be chosen depending on specific problems. In our case, we chose Multi-Quadratic radial basis function.
 $\textit{l} $ represents the Euclidean distance $\|(\boldsymbol{\psi}^{d,n},\boldsymbol{\psi}^{sd,n+1},\boldsymbol{\xi})-(\boldsymbol{\psi}_{j}^{d,n},\boldsymbol{\psi}_{j}^{sd,n+1},\boldsymbol{\xi}_{j})\|$.
$\epsilon$ denotes the shape parameters, which can be optimized using greedy algorithm \cite{xiao2016nonb}.

After the construction of RBF interpolation, we can analytically estimate the gradient at the 'closet' training data points, e.g assuming the $\textit{i}$-th training $(\boldsymbol{\psi}_{i}^{d,n},\boldsymbol{\psi}_{i}^{sd,n+1},\boldsymbol{\xi}_{i})$, by differentiating the RBF as follows,
\begin{align}
\label{eq22}
& \frac{\partial \boldsymbol{\pounds}^{d,n+1}}{\partial \boldsymbol{\xi}}|_{\boldsymbol{\xi}=\boldsymbol{\xi}_{i}} =\sum_{j=1}^{M}\boldsymbol{\omega}_{j}^{d,n+1} \times \notag \\
& \frac{\partial \theta (||(\boldsymbol{\psi}^{d,n},\boldsymbol{\psi}^{sd,n+1},\boldsymbol{\xi})-(\boldsymbol{\psi}_{j}^{d,n},\boldsymbol{\psi}_{j}^{sd,n+1},\boldsymbol{\xi}_{j})||)}{\partial \boldsymbol{\xi}}|_{\boldsymbol{\xi}=\boldsymbol{\xi}_{i}}
\end{align}
\begin{align}
\label{eq23}
& \frac{\partial \boldsymbol{\pounds}^{d,n+1}}{\partial \boldsymbol{\psi}^{sd,n+1}}|_{\boldsymbol{\psi}^{sd,n+1}=\boldsymbol{\psi}_{i}^{sd,n+1}}=\sum_{j=1}^{M}\boldsymbol{\omega}_{j}^{d,n+1}\times \notag \\
& \frac{\partial \theta (||(\boldsymbol{\psi}^{d,n},\boldsymbol{\psi}^{sd,n+1},\boldsymbol{\xi})-(\boldsymbol{\psi}_{j}^{d,n},\boldsymbol{\psi}_{j}^{sd,n+1},\boldsymbol{\xi}_{j})||)}{\partial \boldsymbol{\psi}^{sd,n+1}}|_{\boldsymbol{\psi}^{sd,n+1}=\boldsymbol{\psi}_{i}^{sd,n+1}}
\end{align}
\begin{align}
\label{eq24}
& \frac{\partial \boldsymbol{\pounds}^{d,n+1}}{\partial \boldsymbol{\psi}^{d,n}}|_{\boldsymbol{\psi}^{d,n}=\boldsymbol{\psi}_{i}^{d,n}}=\sum_{j=1}^{M}\boldsymbol{\omega}_{j}^{d,n+1}\times \notag \\
& \frac{\partial \theta (||(\boldsymbol{\psi}^{d,n},\boldsymbol{\psi}^{sd,n+1},\boldsymbol{\xi})-(\boldsymbol{\psi}_{j}^{d,n},\boldsymbol{\psi}_{j}^{sd,n+1},\boldsymbol{\xi}_{j})||)}{\partial \boldsymbol{\psi}^{d,n}}|_{\boldsymbol{\psi}^{d,n}=\boldsymbol{\psi}_{i}^{d,n}}
\end{align}

\begin{table}[h]
\footnotesize
\centering
\caption{Some well-known radial basis functions}\label{tab1}
\begin{tabular}{l l}
\hline
\textbf{Functions} & \textbf{Definition} \\
\hline
Gaussian & $\theta(l)=e^{-(\frac{l}{\epsilon})^{2}}$ \\
Linear Spline & $\theta(l)=l$ \\
Multi-Quadratic & $\theta(l)=\sqrt{l^{2}+\epsilon^{2}}$ \\
Inverse Caddric & $\theta(l)=\frac{1}{l^{2}+\epsilon^{2}}$ \\
Cubic Spline & $\theta(l)=l^{3}$ \\
Thin Plate Spline & $\theta(l)=l^{2}logl$  \\
Inverse Multistory & $\theta(l)=\frac{1}{\sqrt{l^{2}+\epsilon^{2}}}$ \\
\hline
\end{tabular}
\end{table}

Similarly, the well model Eq.\ref{eq16} also can be approximately constructed using RBF interpolation method as follows,
\begin{align}
\label{eq25}
\textbf{y}^{d,m+1} & \approx \boldsymbol{\hbar}^{d.n+1}(\boldsymbol{\psi}^{d,m+1},\boldsymbol{\xi}) \notag \\
& =\sum_{j=1}^{M}\boldsymbol{\varepsilon}_{j}^{d,m+1}\times \theta (\|(\boldsymbol{\psi}^{d,m+1},\boldsymbol{\xi})-(\boldsymbol{\psi}_{j}^{d,m+1},\boldsymbol{\xi}_{j})\|)
\end{align}

And the gradient at the sets of training data points by differentiating the RBF function Eq.\ref{eq25} with respect to $(\boldsymbol{\psi}_{i}^{d,m+1},\boldsymbol{\xi}_{i})$ can be given by
\begin{align}
\label{eq26}
& \frac{\partial \boldsymbol{\hbar}^{d,n+1}}{\partial \boldsymbol{\xi}}|_{\boldsymbol{\xi}=\boldsymbol{\xi}_{i}}=\sum_{j=1}^{M}\boldsymbol{\varepsilon}_{j}^{d,m+1}\times \notag \\
& \frac{\partial \theta (\|(\boldsymbol{\psi}^{d,m+1},\boldsymbol{\xi})-(\boldsymbol{\psi}_{j}^{d,m+1},\boldsymbol{\xi}_{j})\|)}{\partial \boldsymbol{\xi}}|_{\boldsymbol{\xi}=\boldsymbol{\xi}_{i}}
\end{align}
\begin{align}
\label{eq27}
& \frac{\partial \boldsymbol{\hbar}^{d,n+1}}{\partial \boldsymbol{\psi}^{d,m+1}}|_{\boldsymbol{\psi}^{d,m+1}=\boldsymbol{\psi}_{i}^{d,m+1}}=\sum_{j=1}^{M}\boldsymbol{\varepsilon}_{j}^{d,m+1}\times \notag \\
& \frac{\partial \theta (\|(\boldsymbol{\psi}^{d,m+1},\boldsymbol{\xi})-(\boldsymbol{\psi}_{j}^{d,m+1},\boldsymbol{\xi}_{j})\|)}{\partial \boldsymbol{\xi}}|_{\boldsymbol{\psi}^{d,m+1}=\boldsymbol{\psi}_{i}^{d,m+1}}
\end{align}
where, $\boldsymbol{\hbar}^{d,n+1}(\boldsymbol{\psi}^{d,m+1},\boldsymbol{\xi})$ denotes a RBF interpolation function for the simulated measurements $\textbf{y}^{d,m+1}$ at the time level $\textit{m}$+1 for the subdomain $\Omega^{d}$ from the set $(\boldsymbol{\psi}^{d,m+1},\boldsymbol{\xi})$.
$\boldsymbol{\varepsilon}^{d,m+1}$ is a weighting coefficient vector of size $\textit{M}$ (number of training data sets). $\|(\boldsymbol{\psi}^{d,m+1},\boldsymbol{\xi})-(\boldsymbol{\psi}_{j}^{d,m+1},\boldsymbol{\xi}_{j})\|$ is a scalar distance by $\textit{L}_{2} $ norm.
$\theta$ is a set of specific radial basis functions and weighted by a corresponding coefficient $\boldsymbol{\varepsilon}_{j}^{d,m+1}$.

\subsection{Subdomain POD-TPWL algorithm}
By considering the dynamic interaction between neighboring subdomains as in Eq.\ref{eq15}, the coefficients $\boldsymbol{\psi}^{d,n+1}$ for the subdomain $\Omega^{d}$ can be obtained by the modification of Eq.\ref{eq11} as follows,
\begin{align}
\label{eq28}
\boldsymbol{\psi}^{d,n+1} & \approx \boldsymbol{\psi}_{tr}^{d,n+1}+\textbf{E}_{\boldsymbol{\psi}_{tr}}^{d,n+1}(\boldsymbol{\psi}^{d,n}-\boldsymbol{\psi}_{tr}^{d,n}) \notag \\
& +\textbf{E}_{\boldsymbol{\psi}_{tr}}^{sd,n+1}(\boldsymbol{\psi}^{sd,n+1}-\boldsymbol{\psi}_{tr}^{sd,n+1})+\textbf{G}_{\boldsymbol{\xi}_{tr}}^{n+1}(\boldsymbol{\xi}-\boldsymbol{\xi}_{tr})
\end{align}

Coupling domain decomposition and radial basis function interpolation, the derivative matrices required by POD-TPWL for the subdomain $\Omega^{d}$ are estimated as follows
\begin{align}
\label{eq29}
& \textbf{E}_{\boldsymbol{\psi}_{tr}}^{d,n+1} \approx \frac{\partial \boldsymbol{\pounds}^{d,n+1}}{\partial \boldsymbol{\psi}^{d,n}}|_{\boldsymbol{\psi}^{d,n}=\boldsymbol{\psi}_{tr}^{d,n}}, \textbf{E}_{\boldsymbol{\psi}_{tr}}^{sd,n+1} \approx \frac{\partial \pounds^{d,n+1}}{\partial \boldsymbol{\psi}^{sd,n+1}}|_{\boldsymbol{\psi}^{sd,n+1}=\boldsymbol{\psi}_{tr}^{sd,n+1}} \notag \\
& \textbf{G}_{\boldsymbol{\xi}}^{n+1} \approx \frac{\partial \boldsymbol{\pounds}^{d,n+1}}{\partial \boldsymbol{\xi}}|_{\boldsymbol{\xi}=\boldsymbol{\xi}_{tr}}
\end{align}

Similarly, substituting Eq.\ref{eq26}-Eq.\ref{eq27} into Eq.\ref{eq13}, the simulated measurements $\textbf{y}^{d,m+1}$ of the subdomain $\Omega^{d}$ can be reformulated as
\begin{equation}
\label{eq30}
\textbf{y}^{d,m+1} \approx \textbf{y}_{tr}^{d,m+1}+\textbf{A}_{\boldsymbol{\psi}_{tr}}^{d,m+1}(\boldsymbol{\psi}^{d,m+1}-\boldsymbol{\psi}_{tr}^{d,m+1})+\textbf{B}_{\boldsymbol{\xi}_{tr}}^{m+1}(\boldsymbol{\xi}-\boldsymbol{\xi}_{tr})
\end{equation}
\begin{equation}
\label{eq31}
\textbf{A}_{\boldsymbol{\psi}}^{d,m+1} \approx \frac{\partial \boldsymbol{\hbar}^{d,m+1}}{\partial \boldsymbol{\psi}^{d,m+1}}|_{\boldsymbol{\psi}^{d,m+1}=\boldsymbol{\psi}_{tr}^{d,m+1}} , B_{\boldsymbol{\xi}_{tr}}^{m+1} \approx \frac{\partial \boldsymbol{\hbar}^{d,m+1}}{\partial \boldsymbol{\xi}}|_{\boldsymbol{\xi}=\boldsymbol{\xi}_{tr}}
\end{equation}

Our reformulated subdomain POD-TPWL algorithm has three underlying advantages over the traditional POD-TPWL  algorithm: (1) the approximation of the derivative matrices is non-intrusive, e.g, it does not require the modification of legacy code; (2) the implementation of POD-TPWL is local in each subdomain, which has the potential to capture features dominated by local dynamics better than global approximations. Therefore, we could also refer to the subdomain POD-TPWL algorithm as local POD-TPWL; (3) The non-adjacent subdomains almost have no direct dynamic interactions, this kind of subdomain POD-TPWL algorithm can be easily parallelized. Referring to Fig.\ref{fig1}, subdomains $\Omega^{1}, \Omega^{3}, \Omega^{5}, \Omega^{7}, \Omega^{9}$ have no direct interactions, and therefore the subdomain POD-TPWL algorithm can be simultaneously implemented in these five subdomains. This is similar for the subdomains $\Omega^{2}, \Omega^{4}, \Omega^{6}, \Omega^{8}$.

The subdomain POD-TPWL algorithm consists of an offline stage and an online stage. The offline stage describes a computational procedure on how to construct a set of local RBF and estimate the derivative information for each subdomain. Firstly, the solutions of the full-order model are saved as a sequence of snapshots over the whole computational domain and then partitioned into subdomains. From the local sequence of snapshots over each subdomain, a number of local basis vectors is formed using POD, and then unlike the traditional practices that RBF is used to construct a set of surrogates for each subdomain, we use RBF to estimate the derivative matrices for each subdomain. Finally, those estimated derivative matrices are substituted into POD-TPWL algorithm to form a reduced-order linear model in each subdomain. While the online stage describes how to iteratively implement the subdomain POD-TPWL where the dynamic interactions between a subdomain and its surrounding subdomains are considered. Referring to Eq.\ref{eq28}, the variables of one subdomain at current time level can be linearized around the variables of this subdomain at previous timestep and variables of neighboring subdomains at current timestep, which have not been determined. Thus, some additional iterative steps are needed.

\subsection{Sampling Strategy}
In our proposed subdomain POD-TPWL algorithm, training points are required for both RBF interpolation and to construct the POD basis. For POD, the snapshot matrix generated from the training simulations are expected to sufficiently preserve the dynamic behavior of the system. For RBF interpolation, the training points are selected to compute the derivative matrices. The procedure on how to choose these training points will be described here.

$\textsl{Sampling strategy for POD}$. A set of parameters is initially sampled and used as input for full-order model (FOM) simulations from which a snapshot matrix is constructed. The singular value spectrum is computed for this initial set of samples. The number of samples is then increased one at a time, i.e. adding one FOM simulation, and the SVD is recomputed, until no significant changes are observed in the singular value spectrum.

$\textsl{Sampling strategy for RBF}$. The accuracy of the RBF interpolation will be reduced if too few data points are chosen, while the computational cost increases with the number of data points, which will be prohibitive if too many points are chosen. To limit the number of FOM simulations used to construct the interpolation model for the POD coefficients we use 2-sided perturbation of each coefficient $\xi_{j}$ resulting in $2 \times l_\beta + 1$ points. In some experiments we add additional points by simultaneous random sampling of perturbations $\Delta \boldsymbol{\xi}$. An alternative could be use to use Smolyak sparse grid sampling \cite{smolyak1963quadrature}.

\subsection{Adjoint-based history matching algorithm}
After linearizing the original full-order model to a reduced-order linear model, because of the linearity in the reduced-order space, the implementation of the adjoint model is easily realized. It is convenient to incorporate this reduced-order linear model established using the subdomain POD-TPWL into the adjoint-based reservoir history matching.

The cost function in the reduced-order space can be given by reformulating the Eq.\ref{eq3} as follows,
\begin{align}
\label{eq33}
& \jmath(\boldsymbol{\xi}) = \frac{1}{2}(\boldsymbol{\beta}_{b}+\boldsymbol{\phi}_{\boldsymbol{\beta}} \boldsymbol{\xi}-\boldsymbol{\beta}_{p})^{T}{\textbf{R}_{p}}^{-1}(\boldsymbol{\beta}_{b}+\boldsymbol{\phi}_{\boldsymbol{\beta}} \boldsymbol{\xi}-\boldsymbol{\beta}_{p})  \notag \\
& +\frac{1}{2}\sum_{d=1}^{S}\sum_{m=1}^{N_{0}}[\textbf{d}_{o}^{d,m}-\textbf{y}_{tr}^{d,m}-\textbf{A}_{\boldsymbol{\psi}_{tr}}^{d,m}(\boldsymbol{\psi}^{d,m}-\boldsymbol{\psi}_{tr}^{d,m})-\textbf{B}_{\boldsymbol{\xi}_{tr}}^{m}(\boldsymbol{\xi}-\boldsymbol{\xi}_{tr})]^{T}   \notag \\
& {\textbf{R}_{m}}^{-1}[\textbf{d}_{o}^{d,m}-\textbf{y}_{tr}^{d,m}-\textbf{A}_{\boldsymbol{\psi}_{tr}}^{d,m}(\boldsymbol{\psi}^{d,m}-\boldsymbol{\psi}_{tr}^{d,m})-\textbf{B}_{\boldsymbol{\xi}_{tr}}^{m}(\boldsymbol{\xi}-\boldsymbol{\xi}_{tr})]
\end{align}
and its gradient is
\begin{align}
\label{eq34}
& [\frac{d\jmath}{d\boldsymbol{\xi}}]^{T} =(\boldsymbol{\phi}_{\boldsymbol{\beta}} )^{T}
{\textbf{R}_{p}}^{-1}(\boldsymbol{\beta}_{b}+\boldsymbol{\phi}_{\boldsymbol{\beta}} \boldsymbol{\xi}-\boldsymbol{\beta}_{p}) \notag \\
& -\sum_{d=1}^{S}\sum_{m=1}^{N_{0}}[\textbf{B}_{\boldsymbol{\xi}_{tr}}^{m}]^{T}
{\textbf{R}_{m}}^{-1}[\textbf{d}_{o}^{d,m}-\textbf{y}_{tr}^{d,m}-\textbf{A}_{\boldsymbol{\psi}_{tr}}^{d,m}(\boldsymbol{\psi}^{d,m}-\boldsymbol{\psi}_{tr}^{d,m}) \notag \\
& -\textbf{B}_{\boldsymbol{\xi}_{tr}}^{m}(\boldsymbol{\xi}-\boldsymbol{\xi}_{tr})]-\sum_{d=1}^{S}\sum_{n=1}^{N}[\textbf{G}_{\boldsymbol{\xi}_{tr}}^{n}]^{T}\boldsymbol{\lambda}^{d,n}
\end{align}
where $\boldsymbol{\lambda}^{d,n}$ is obtained as the solution of the adjoint model for the subdomain $\Omega^{d}$ is given by
\begin{align}
\label{eq35}
& [\textbf{I}-(\textbf{E}_{\boldsymbol{\psi}_{tr}}^{d,n})^{T}]\boldsymbol{\lambda}^{d,n}  = \sum_{d=1}^{S}[\textbf{A}_{\boldsymbol{\psi}_{tr}}^{d,n}]^{T}
{\textbf{R}_{n}}^{-1}[\textbf{d}_{o}^{d,n}-\textbf{y}_{tr}^{d,n} \notag \\
& -\textbf{A}_{\boldsymbol{\psi}_{tr}}^{d,n}(\boldsymbol{\psi}^{d,n}-\boldsymbol{\psi}_{tr}^{d,n})-\textbf{B}_{\boldsymbol{\xi}_{tr}}^{n}(\boldsymbol{\xi}-\boldsymbol{\xi}_{tr})]+[\textbf{E}_{\boldsymbol{\psi}_{tr}}^{sd,n}]^{T}\boldsymbol{\lambda}^{d,n+1}
\end{align}
The minimization of the cost function (Eq.\ref{eq33}) can be performed using a steepest descent algorithm \cite{Nocedal1999Numerical} and is stopped when either one of the following stopping criteria is satisfied
\begin{itemize}
\item No more change in the cost function,

\begin{equation}
\label{eq36}
\frac{|\jmath(\boldsymbol{\xi}^{k+1})-\jmath(\boldsymbol{\xi}^{k})|}{max{|\jmath(\boldsymbol{\xi}^{k+1})|,1}}< \eta_{\jmath}
\end{equation}

\item No more change in the estimate of parameters,

\begin{equation}
\label{eq37}
\frac{|\boldsymbol{\xi}^{k+1}-\boldsymbol{\xi}^{k}|}{max{|\boldsymbol{\xi}^{k+1}|,1}}< \eta_{\boldsymbol{\xi}}
\end{equation}

\item The maximum number of iterations has been reached. i.e

\begin{equation}
\label{eq38}
k <= N_{max}
\end{equation}

\end{itemize}
where $\eta_{\jmath}$ and $\eta_{\boldsymbol{\xi}}$ are predefined error constraints and $N_{max}$ is the maximum number of iterations.

As mentioned in \cite{kaleta2011model}, the solution of the reduced and linearized minimization problem based on Eq.\ref{eq33} is not necessarily the solution of the original problem based on Eq.\ref{eq3}. Therefore an additional stopping criterion should be introduced for the original model as follows \cite{Tarantola2005Inverse},
\begin{equation}
\label{eq39}
N_{d}N_{0}-2\sqrt{2N_{d}N_{0}} \leqslant 2J(\boldsymbol{\beta}^k) \leqslant N_{d}N_{0}+2\sqrt{2N_{d}N_{0}}
\end{equation}
where, $\textit{N}_{0}$ is the number of timesteps where the measurements are taken, $ N_{d} $ is the number of measurements at each timestep, $\boldsymbol{\beta}^{k}$ represents the
updated parameters vector at the $\textit{k}$-th outer-loop. $\textit{J}$ is the cost function computed as Eq.\ref{eq3}.

\begin{figure}[h]
\centering\includegraphics[width=1\linewidth]{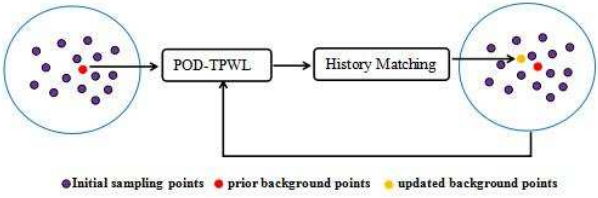}
\caption{The illustration of the reconstruction method of subdomain POD-TPWL algorithm}\label{fig2}
\end{figure}

\begin{figure}[!h]
\centering
\includegraphics[width=1\linewidth]{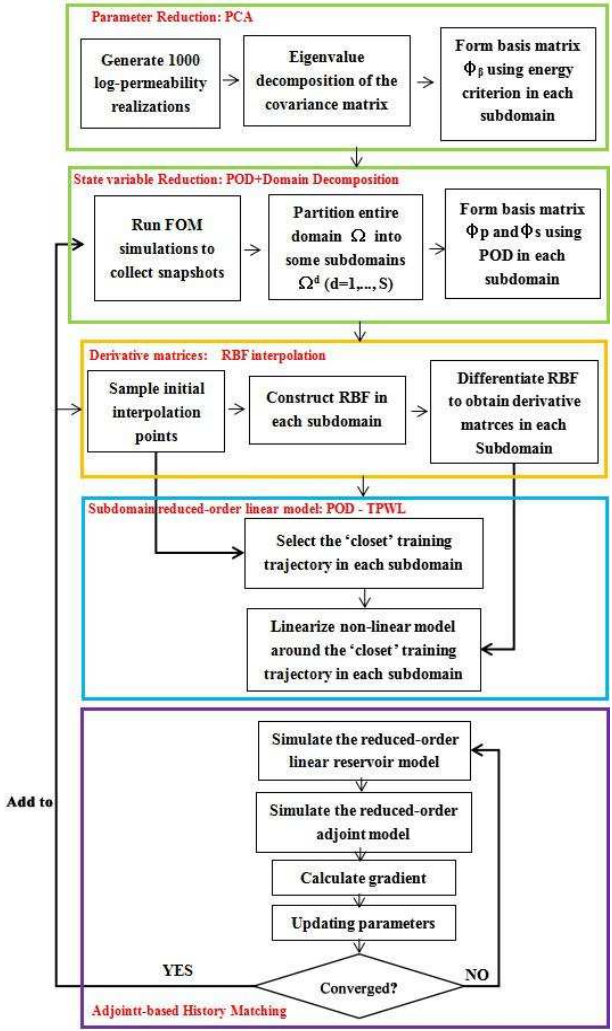}
\caption{The illustration of adjoint-based history matching using subdomain POD-TPWL algorithm}\label{fig3}
\end{figure}

If the objective function does not obey the stopping criterion as Eq.\ref{eq39}, then additional outer-loops are required to reconstruct new reduced-order linear models using the updated parameters. Since the dynamic patterns are mainly dominated by the model parameterization, whenever the parameters are changed, a new background state and a new set of patterns need to be identified. Thereafter, a new reduced-order linear model is built and the aforementioned iterative inner-loop is performed again.

Our proposed non-intrusive subdomain POD-TPWL has computational advantages over the traditional construction of reduced-order linear models using perturbation-based finite-difference method proposed in \cite{kaleta2011model}, especially when the reduced-order linear model is required to be reconstructed for each outer-loop.
Instead of re-perturbing the parameter and state variables one by one to approximate the derivative matrices as proposed in \cite{kaleta2011model}, which would require an additional  ($l_{p}$+$l_{s}$+$l_{\boldsymbol{\beta}}$+1) full order
model (FOM) simulations, our algorithm runs only one additional FOM simulation using updated parameters. The updated parameters and simulated snapshots are added into the previous group of sampling interpolation points and corresponding snapshots. The derivative matrices for the updated parameters are approximated based on the updated group of interpolation points and snapshots. The overall workflow has been summarized conceptually in Fig.\ref{fig2}. The individual steps of the history matching algorithm described in this section are summarized in the flow chart presented in Fig.\ref{fig3}.

\section{Numerical experiments and Discussion}
In this section, some numerical experiments are presented that aim to demonstrate and evaluate our proposed adjoint-based history matching algorithm. In our numerical experiments, MRST, a free open-source software for reservoir modeling and simulation\cite{lie2012open}, is used to run the full-order model simulations.

\subsection{Description of model settings}
\begin{figure}[!h]
\centering
\includegraphics[width=0.8\linewidth]{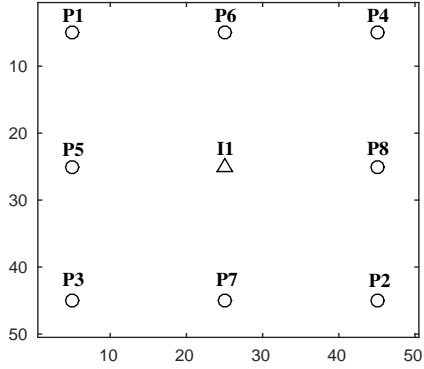}
\caption{The well placement in the 2-D reservoir model for case 1}\label{fig4}
\end{figure}

A 2D heterogeneous oil-water reservoir is considered with two-phase imcompressible flow dynamics. The reservoir contains 8 producers and 1 injector, which are labeled as $P_{1}$ to $P_{8}$, and $I_{1}$ respectively, see Fig.\ref{fig4}. Detailed information about the reservoir geometry, rock properties, fluid properties, and well controls are summarized in Table \ref{tab2}.

\begin{table}[h]
\footnotesize
\centering
\caption{Experiment settings using MRST fro case 1}\label{tab2}
\begin{spacing}{1.25}
\begin{tabular}{l l}
\hline
\textbf{Description} & \textbf{Value} \\
\hline
Dimensions & 50 $\times$ 50 $\times $1 \\
Grid cell size & 20 $\times$ 20 $\times$ 10 \\
Number of wells & 8 producers, 1 injector \\
Fluid density & 1014 kg/$m^{3}$, 859 kg/$m^{3}$ \\
Fluid viscosity & 0.4 mP$\cdot$s, 2 mP$\cdot$s \\
Initial pressure & 30 MPa  \\
Initial saturation & S$_{o}$=0.80,  S$_{w}$=0.20 \\
Connate water saturation & $S_{wc}$=0.20 \\
Residual oil saturation & $S_{or}$=0.20  \\
Corey exponent, oil & 4.0   \\
Corey exponent, water & 4.0  \\
Injection rate & 200$m^{3}$/d  \\
BHP & 25MPa  \\
History production time & 5 year  \\
Prediction time & 10 year  \\
Timestep & 0.1 year  \\
Measurement timestep & 0.2 year  \\
\hline
\end{tabular}
\end{spacing}
\end{table}

\subsection{Reduced model construction}
We generate an ensemble of 1000 Gaussian-distributed realizations of log-permeability. We also assume that the generated log-permeability fields are not conditioned to
the permeability values at the well locations. The log-permeability fields and the corresponding porosity fields are described by the following statistics:
\begin{equation}
\label{eq40}
\sigma_{\boldsymbol{\beta}}=5
\end{equation}
\begin{equation}
\label{eq41}
\textbf{C}_{\boldsymbol{\beta}}(x_{i1,j1};y_{i2,j2}]=\sigma_{\boldsymbol{\beta}}^{2}e^{-[(\frac{|x_{i1}-x_{i2}|}{\chi_{x}})^{2}+(\frac{|y_{i1}-y_{i2}|}{\chi_{y}})^{2}]}
\end{equation}
\begin{equation}
\label{eq42}
\frac{\chi_{x}}{L_{x}}=0.2, \frac{\chi_{y}}{L_{y}}=0.2
\end{equation}
\begin{equation}
\label{eq43}
\boldsymbol{\phi}=0.25(\frac{e^{\boldsymbol{\beta}}}{200})^{0.1}
\end{equation}

\begin{figure}[!h]
\centering
\subfloat['True' model in original full-order space]%
  {\includegraphics[width=0.7\linewidth]{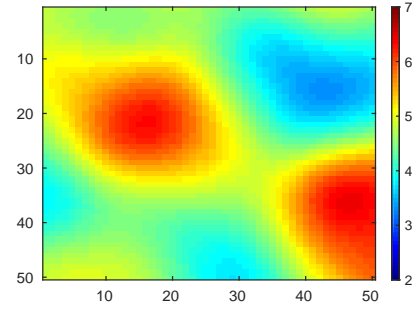}}\
\subfloat[Projected 'true' model in reduced-order subspace]%
  {\includegraphics[width=0.7\linewidth]{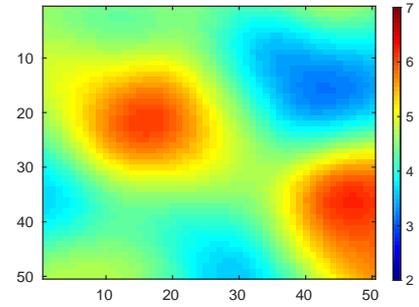}}\
\caption{Comparison of the 'true' reservoir model in full-order space and in reduced-order space for Case 1).}\label{fig5}
\end{figure}

\begin{figure}[!h]
\centering\includegraphics[width=1\linewidth]{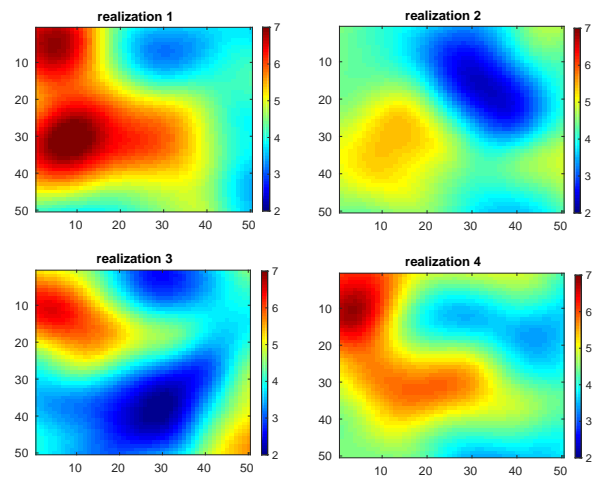}\\
\caption{Examples of realizations of the log-permeability field generated from PCA coefficients sampled randomly from the set ${-1, 1}$ for Case 1.}\label{fig6}
\end{figure}

Here, $\sigma_{\boldsymbol{\beta}} $ is the standard deviation of log-permeability $\boldsymbol{\beta}$; $\textit{C}_{\boldsymbol{\beta}}$ is the covariance of $\boldsymbol{\beta}$; $x_{i1,j1}$=($x_{i1}$,$y_{j1}$)
denotes the coordinates of a grid block; $\chi_{x}$ (or $\chi_{y}$) is the correlation length in $\textit{x}$ (or $\textit{y}$) direction;
and $\textit{L}_{x}$ (or $\textit{L}_{y}$) is the domain length in $\textit{x}$ (or $\textit{y}$) direction.
The background log-permeability $\boldsymbol{\beta}_{b}$ is taken as the average of the 1000 realizations.
One of the realizations was considered to be the truth, and is illustrated in Fig.\ref{fig5}(a).
The permeability field was parameterized using KL-expansion and about 95\% energy is maintained,
resulting in 18 permeability patterns with $l_{\boldsymbol{\beta}} = 18$ corresponding independent PCA coefficients,
which are used in the workflow as a low-dimensional representation of the 2500 grid block permeability values.
Fig.\ref{fig5}(b) shows the projection of the 'true' permeability field in this low-dimensional subspace which shows that the truth can be almost perfectly reconstructed in this subspace.
Four realizations for log-permeability field generated are additionally shown in Fig.\ref{fig6}.

After having reduced the parameter space, the next step is to reduce the reservoir model. The first step is to generate a set of training runs from which snapshots will be taken.
Since the required number of training runs is not known a priori we follow the following procedure: (1) generate a sample PCA coefficient vector by sampling from the set ${-1,1}$,
(2) run a full-order model simulation with these parameters, (3) extract snapshots and add to a snapshot matrix,
(4) compute the singular value decomposition of the snapshot matrix (5) repeat steps (1) to 4) until changes in the singular values are insignificant.
For this case this produced a set of 15 training runs and 240 snapshots for pressure and saturation each.

\begin{figure}[!h]
\centering
\includegraphics[width=0.7\linewidth]{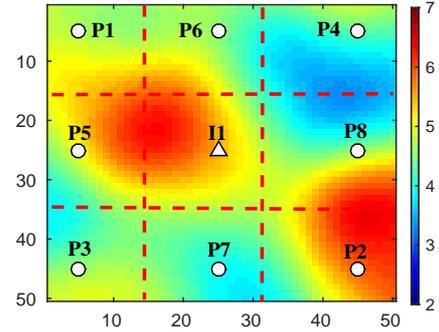}
\caption{Illustration of the applied domain decomposition for Case 1.}\label{fig7}
\end{figure}

The next step is to build a local RBF interpolation model. We divide the entire domain into 9 rectangle subdomains as illustrated in Fig.\ref{fig7}.
The choice of subdomains is fairly arbitrary at this point since we have no formal algorithm to determine the best number and design of the subdomains.
The previously collected global snapshots for pressures and saturations are divided into local snapshots.
For each subdomain, two separate eigenvalue problems for pressure and saturation are solved using POD.
The number of the reduced parameters and state patterns for each subdomain and the global domain are listed in Table \ref{tab3} where 95\% and 90\% of energy are preserved for the pressure and saturation respectively in each subdomain.

After implementing the KL-expansion, original parameters are represented by a set of independent Gaussian random variables with zero mean and unit variance.
In our case, the initial 37 sampling points are selected within interval [-1,1] as described in subsection of sampling strategy,
where the \textit{j}-th element $\boldsymbol{\xi}_{i}^{j}$ of the \textit{i}-th PCA coefficient vector $\boldsymbol{\xi}^{j}$
is perturbed sequentially in 2 opposite directions (positive and negative) by a specific amplitude perturbation $\Delta$ $\boldsymbol{\xi}_{i}^{j}$.

\begin{figure}[!h]
\centering\includegraphics[width=1\linewidth]{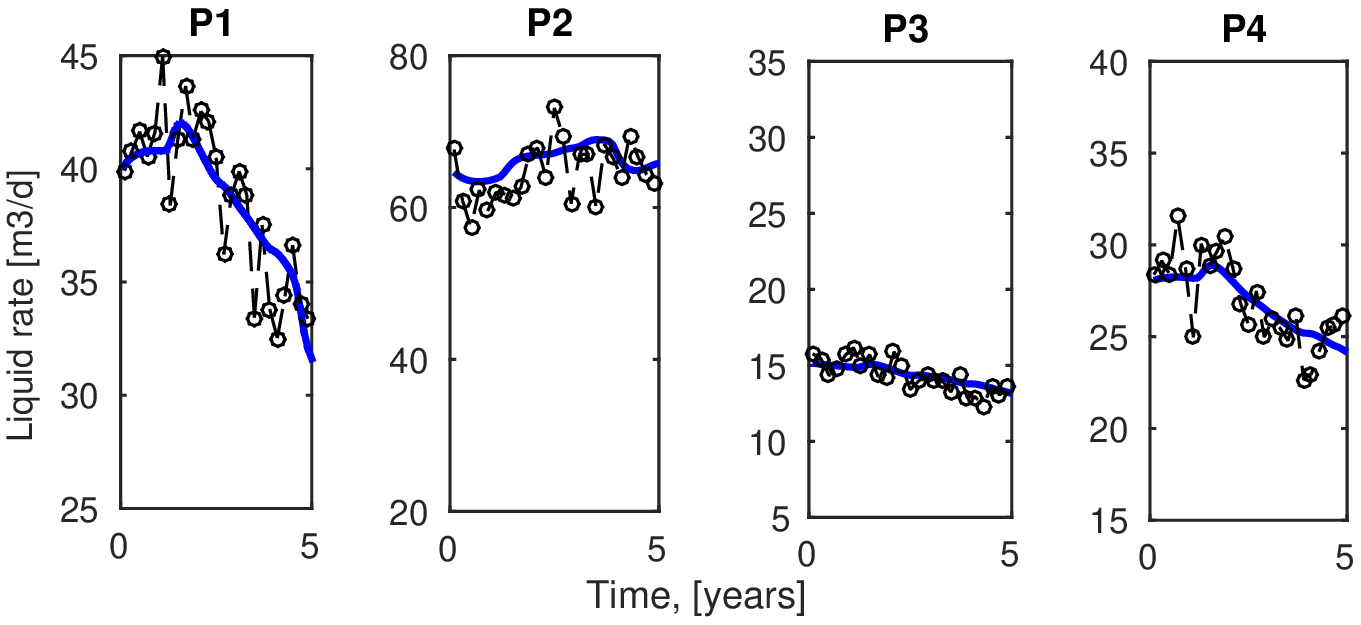}\\
\centering\includegraphics[width=1\linewidth]{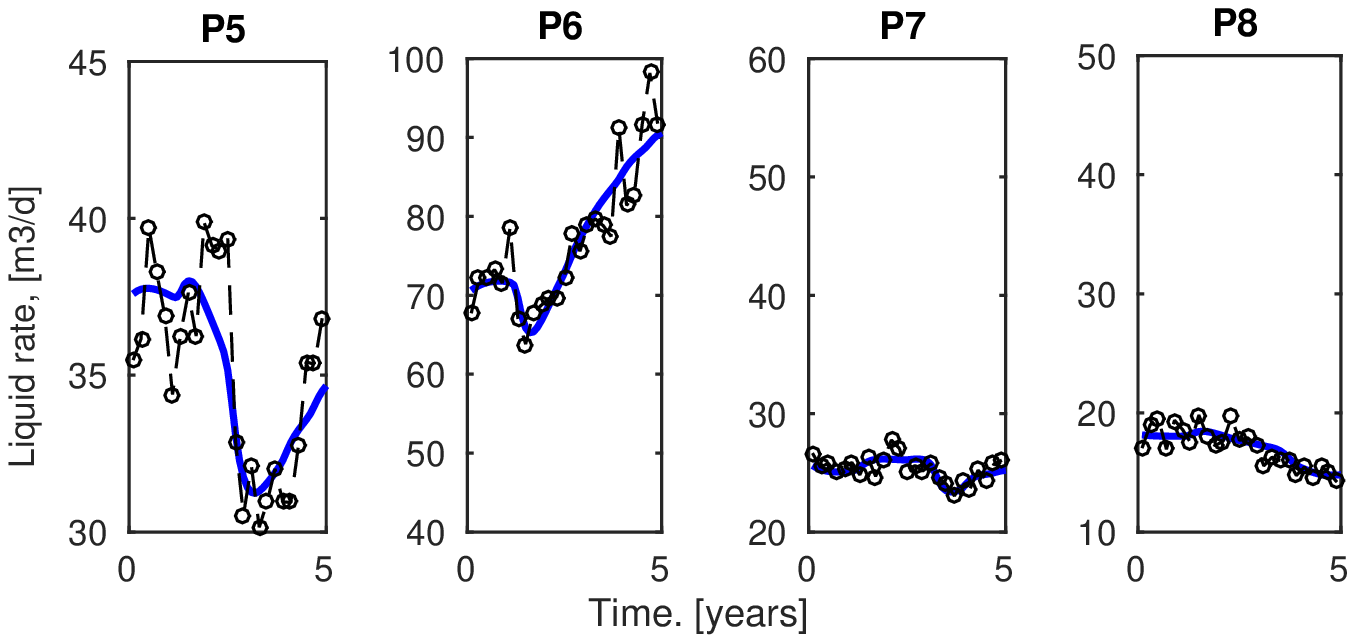}\\
\centering\includegraphics[width=1\linewidth]{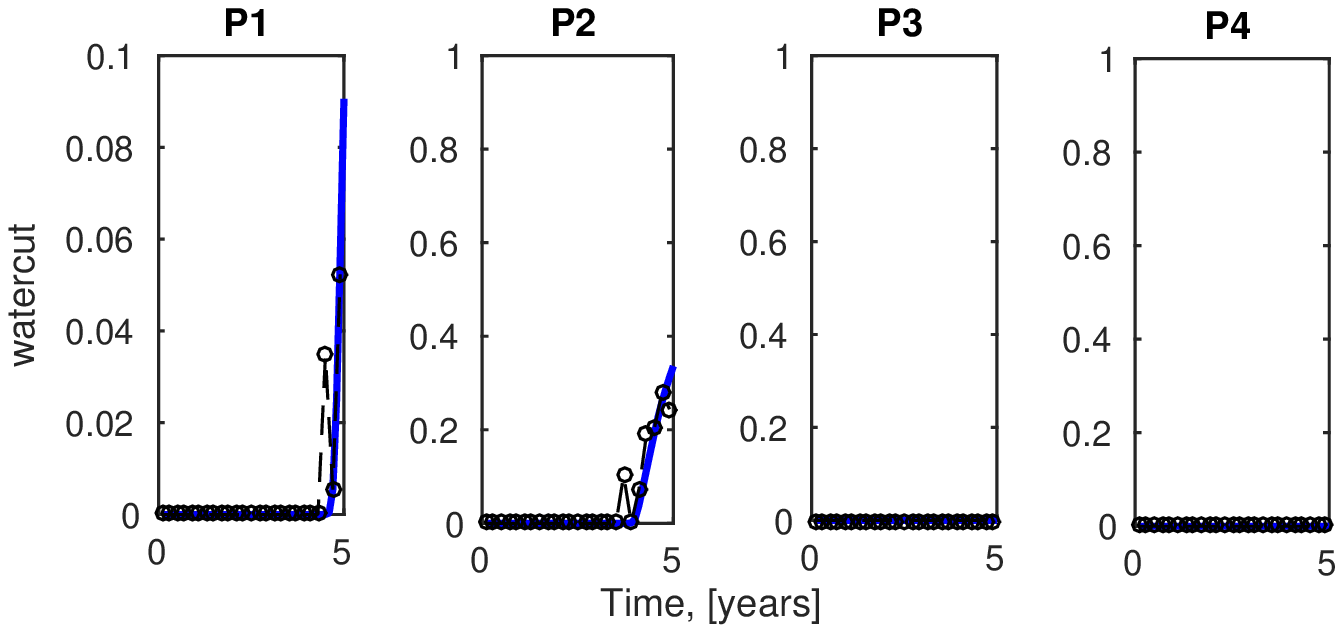}\\
\centering\includegraphics[width=1\linewidth]{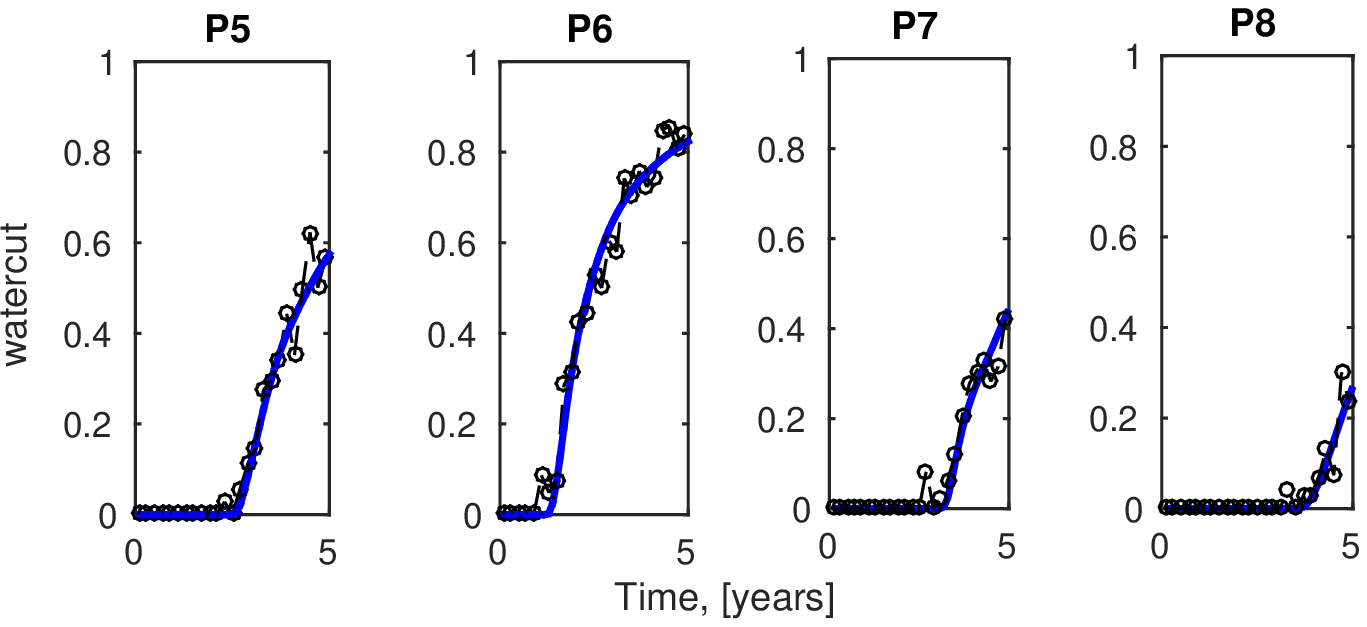} \\
\centering\includegraphics[width=1\linewidth]{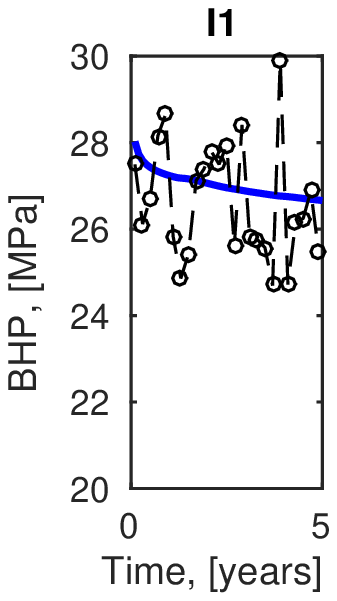} \\
\caption{Measured quantities for case 1. Blue solid line: reference model (truth). Black dashed line: noisy data.}\label{fig8}
\end{figure}

The history period is 5 years during which observations are taken from 8 producers and 1 injector every two model timesteps (nearly 73 days) resulting in 25 time instances. Noisy observations
are generated from the model with the "true" permeability field and include bottom-hole
pressures (BHP) in the injector and fluid rates and water-cut (WCT) in the producers. As a result we have 200 fluid rates and 200 WCT values measured in the
producers and 25 bottom-hole pressures measured in the injector, which gives in total 425 measurement data.
Normal distributed independent measurement noise with a standard deviation equal to 5\% of the 'true' data value, was added to all observations. The generated measurements are shown in Fig.\ref{fig8}.

To analyze the results, we define two error measures based on data misfits $e_{obs}$ and parameter misfits $e_{\boldsymbol{\beta}}$ as follows,
\begin{equation}
\label{eq44}
e_{obs}=\sqrt{\frac{\sum_{i=1}^{N_{o}}\sum_{j=1}^{N_{d}}(\textbf{d}_{obs}^{i,j}-\textbf{d}_{upt}^{i,j})^{2}}{N_{o}N_{d}}}
\end{equation}
\begin{equation}
\label{eq45}
e_{\boldsymbol{\beta}}=\sqrt{\frac{\sum_{i=1}^{N_{g}}(\boldsymbol{\beta}_{true}^{i}-\boldsymbol{\beta}_{upt}^{i})^2}{N_{g}}}
\end{equation}
where, $\textbf{d}_{obs}^{i,j}$ and $\textbf{d}_{upt}^{i,j}$ represent the measurements and simulated data using the updated model respectively; $\boldsymbol{\beta}_{true}^{i}$ and $\boldsymbol{\beta}_{upt}^{i}$ denote the
grid block log-permeability from the 'true' model and updated model respectively.

\begin{table}[!h]
\footnotesize
\centering
\caption{Summary of the number of reduced variables for the global domain and after domain decomposition for case 1 (Note:\textbf{s} refers to saturation, \textbf{p} refers to pressure)}\label{tab3}
\begin{tabular}{|c|c|c|c|c|c|c|c}
\hline
\multicolumn{4}{|c|}{Domain Decomposition} & \multicolumn{3}{|c|}{Global Domain} \\
\hline
SD & $\boldsymbol{\beta}$ & \textbf{s} & \textbf{p} & $\boldsymbol{\beta}$ & \textbf{s} & \textbf{p} \\
\hline
1 & \multirow{9}{*}{18} & 14 & 7  &\multirow{9}{*}{18} & \multirow{9}{*}{72} &\multirow{9}{*}{36}\\
2 & \multirow{9}{*}{} & 13 & 6  &\multirow{9}{*}{} & \multirow{9}{*}{} &\multirow{9}{*}{}\\
3 & \multirow{9}{*}{} & 12 & 5  &\multirow{9}{*}{} & \multirow{9}{*}{} &\multirow{9}{*}{}\\
4 & \multirow{9}{*}{} & 13 & 4  &\multirow{9}{*}{} & \multirow{9}{*}{} &\multirow{9}{*}{}\\
5 & \multirow{9}{*}{} & 16 & 7  &\multirow{9}{*}{} & \multirow{9}{*}{} &\multirow{9}{*}{}\\
6 & \multirow{9}{*}{} & 14 & 6  &\multirow{9}{*}{} & \multirow{9}{*}{} &\multirow{9}{*}{}\\
7 & \multirow{9}{*}{} & 13 & 5  &\multirow{9}{*}{} & \multirow{9}{*}{} &\multirow{9}{*}{}\\
8 & \multirow{9}{*}{} & 15 & 6  &\multirow{9}{*}{} & \multirow{9}{*}{} &\multirow{9}{*}{}\\
9 & \multirow{9}{*}{} & 12 & 5  &\multirow{9}{*}{} & \multirow{9}{*}{} &\multirow{9}{*}{}\\
\hline
Total & 18 & 122 & 51 & 18 & 72 & 36 \\
\hline
\end{tabular}
\end{table}

\subsection{History matching results}
\begin{figure}[!h]
\centering\includegraphics[width=1\linewidth]{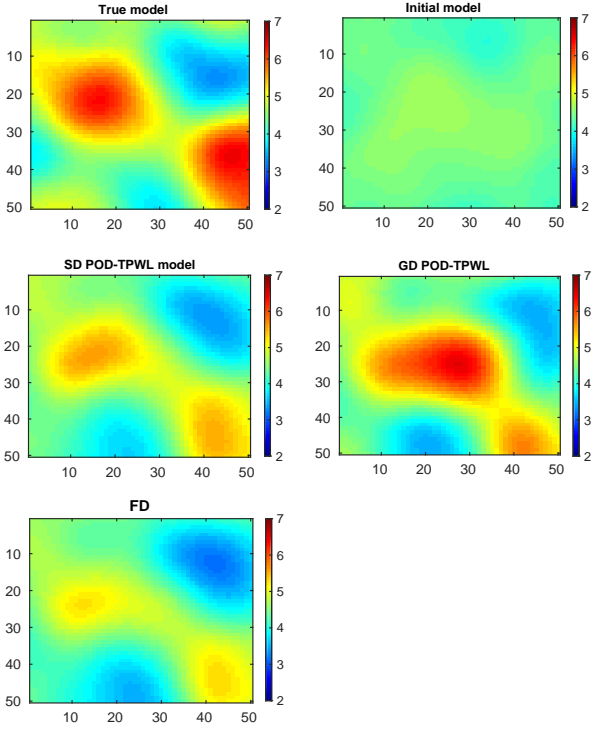}\\
\caption{True, initial and updated log-permeability fields using SD POD-TPWL, GD POD-TPWL, and the FD method for case 1.}\label{fig9}
\end{figure}

\begin{figure}[!h]
\centering
\subfloat[FD method]{%
  \includegraphics[width=0.8\linewidth]{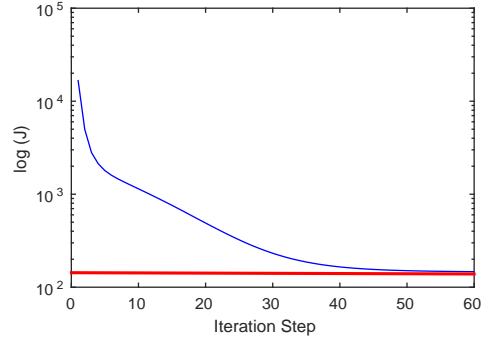}}\hfill
\subfloat[subdomain POD-TPWL]{%
  \includegraphics[width=0.8\linewidth]{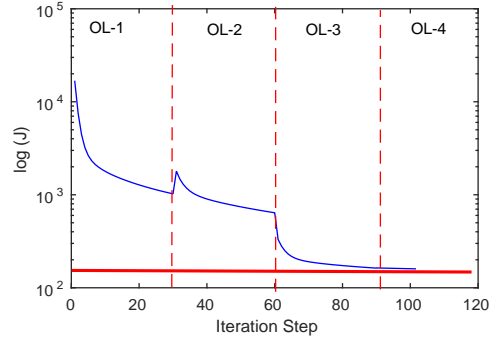}}\\
\caption{Cost function value decrease using (a) finite-difference method, and (b) subdomain POD-TPWL for case 1. OL-i means the i-th outer-loop.}\label{fig10}
\end{figure}

\begin{figure}[!h]
\centering\includegraphics[width=1\linewidth]{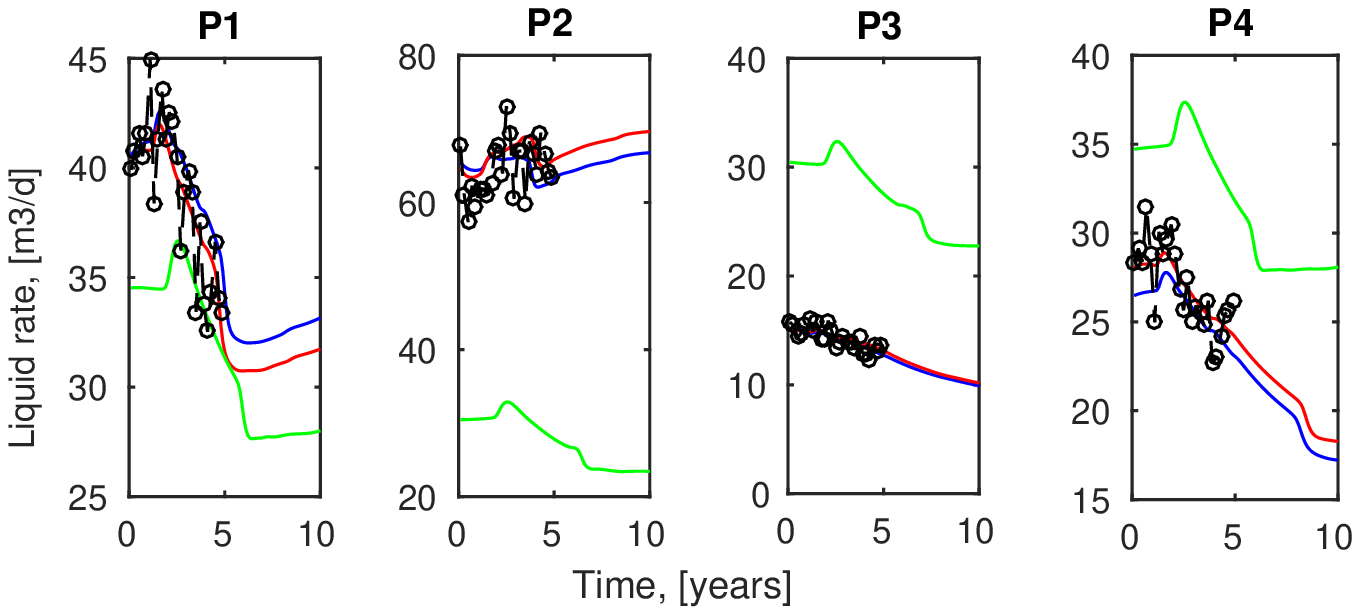}\\
\centering\includegraphics[width=1\linewidth]{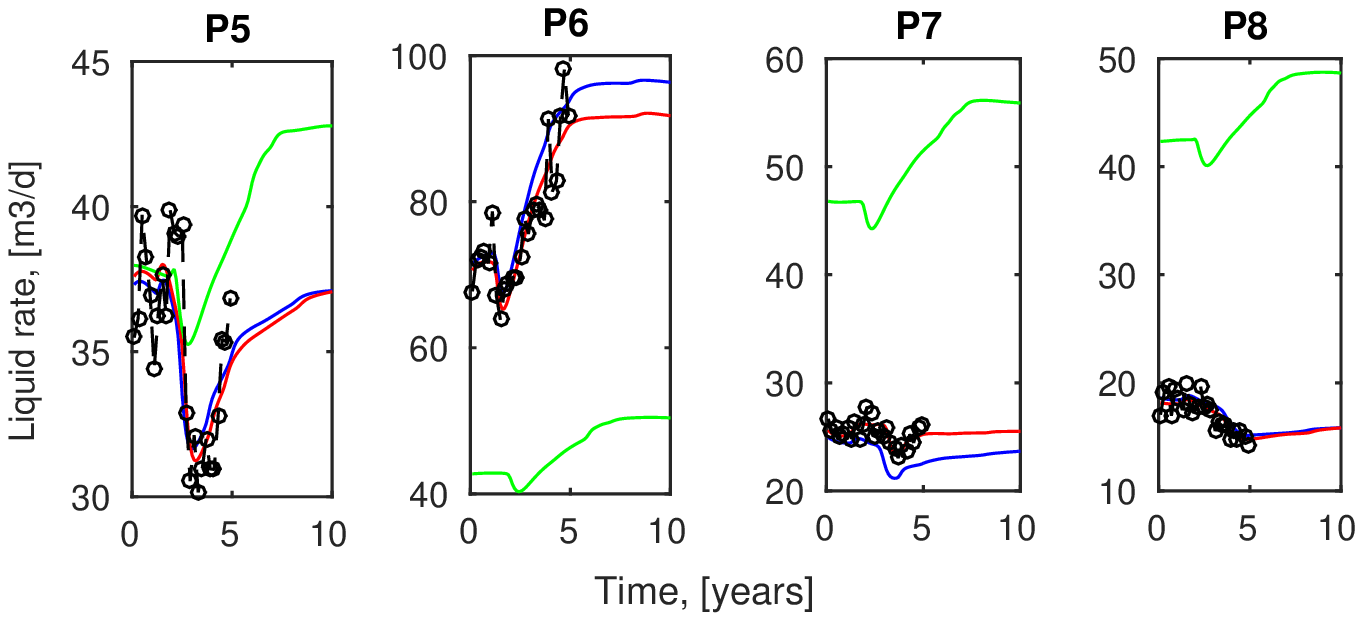}\\
\centering\includegraphics[width=1\linewidth]{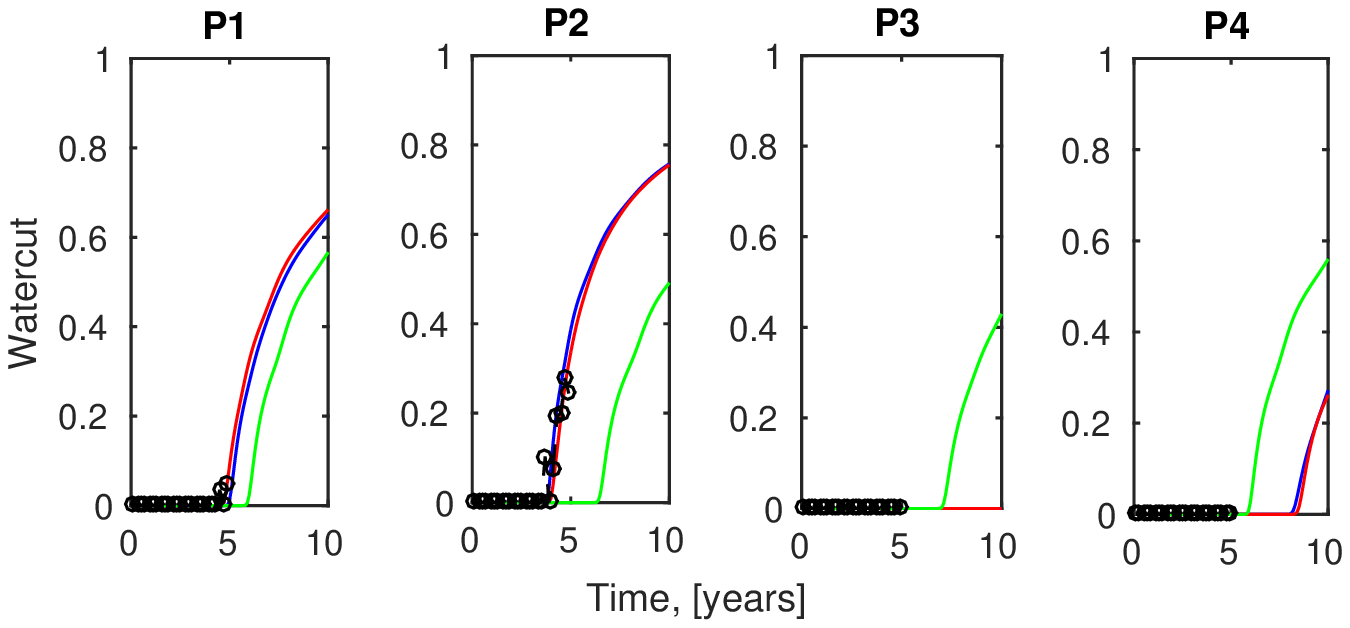}\\
\centering\includegraphics[width=1\linewidth]{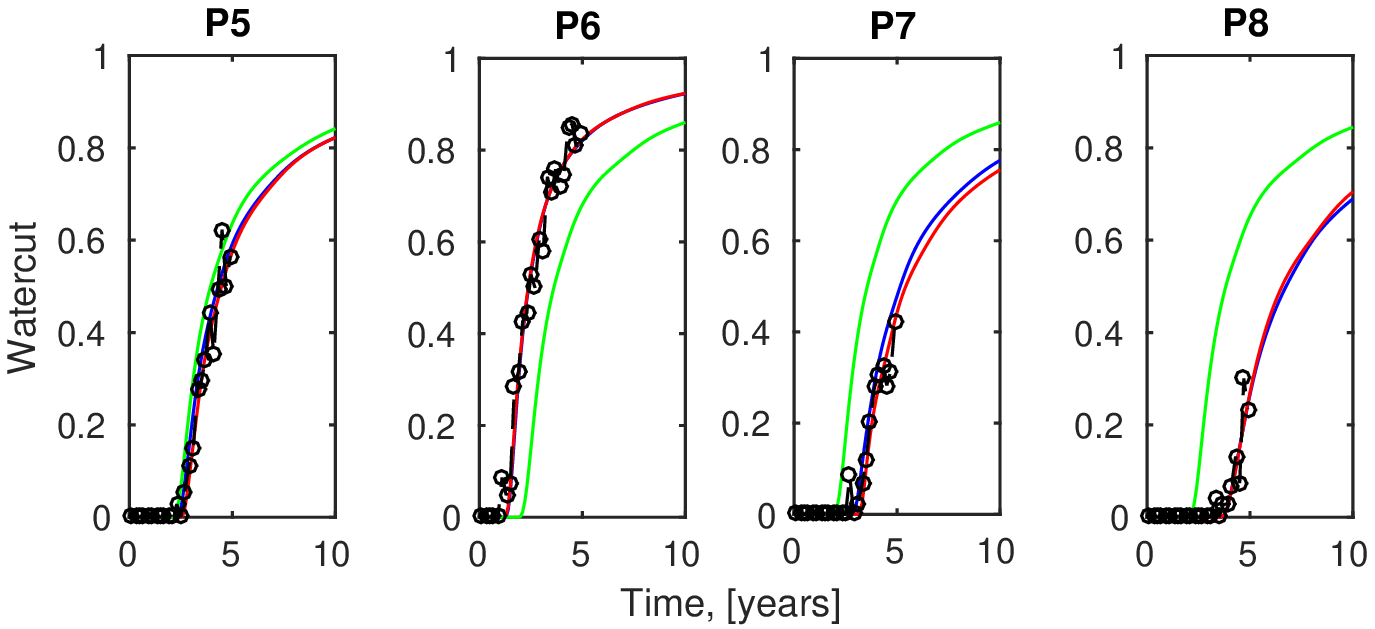}\\
\caption{Forecasts of the liquid rate and water-cut for case 1 from the initial model (green line), the reference model (blue line), and the model updated using the SD POD-TPWL (red line). Measured data are indicate by open circles.}\label{fig11}
\end{figure}

Figures \ref{fig9}, \ref{fig10} and \ref{fig11} and Table \ref{tab4} show the results of the first numerical experiments, including the updated log-permeability field,
the value of cost function at each iteration and the mismatch between observed data and predictions. To demonstrate the performance of our proposed methodology, we compared the
results with those of finite-difference (FD) based history matching algorithm without domain decomposition and model order reduction.
The total computational cost of any minimization problem strongly depends on the number of parameters. In our work, for a fair comparison,
we use the same reparameterization to reduce the number of parameters and implement the finite-difference based history matching in this reduced-order parameter subspace.
The cost function for finite-difference based history matching can be defined as follows,
\begin{align}
\label{eq46}
& J(\boldsymbol{\xi}) = \frac{1}{2}(\boldsymbol{\beta}_{b}+\boldsymbol{\phi}_{\boldsymbol{\beta}} \boldsymbol{\xi}-\boldsymbol{\beta}_{p})^{T}{\textbf{R}_{p}}^{-1}(\boldsymbol{\beta}_{b}+\boldsymbol{\phi}_{\boldsymbol{\beta}} \boldsymbol{\xi}-\boldsymbol{\beta}_{p})  \notag \\
& +\frac{1}{2}\sum_{m=1}^{N_{0}}(\textbf{d}_{o}^{m}-\textbf{h}^{m}(\textbf{x}^{m},\boldsymbol{\xi}))^{T}{\textbf{R}_{m}}^{-1}(\textbf{d}_{o}^{m}-\textbf{h}^{m}(\textbf{x}^{m},\boldsymbol{\xi}))
\end{align}

The finite-difference method is used to compute the numerical gradient of the cost function as Eq.\ref{eq46} with respect to 18 PCA coefficients.
A FD gradient is determined by one-sided perturbation of each of the 18 PCA coefficients. Thus, 19 full-order model (FOM) simulations are required for each iteration step.
The stopping criteria are set $\eta_{\jmath} =10^{-4} $, $\eta_{\boldsymbol{\xi}} =10^{-3} $, and $N_{max}$=30. As can be seen from Fig.\ref{fig10} and Table \ref{tab4}, the model-reduced approach needs 55 full order
model (FOM simulations, consiting of 15 FOM simulations to collect the snapshots and 37 FOM simulations to construct the initial reduced-order linear model in the first outer-loop.
The remaining 3 FOM simulations are used to reconstruct the reduced-order linear models in the next three outer-loops and to calculate the cost function as Eq.\ref{eq3} in the original space.
Fig.\ref{fig9} shows the true, initial and final estimates of log-permeability field. In this case, the main geological structures of the the 'true' model can be reconstructed with both methods.
However, the parameter estimates obtained with proposed methodology more accurately reproduce the true amplitudes than
those obtained with the classic finite-difference based history matching. From Fig.\ref{fig11} and Table \ref{tab4},
we can both qualitatively and quantitatively observe that the history matching process results in an improved prediction in all of the eight production wells.
Fig.\ref{fig9} illustrates the data match of fluid rate and water-cut up to year 5 and an additional 5-year prediction until year 10 for all 8 producers.
The prediction based on the initial model is far from that of the true model.
After history matching, the predictions from the updated model match the observations very well.
Also the prediction of water breakthrough time is imporved for all of the production wells, also for wells that show water breakthrough only after the history period.


\begin{table}[!h]
\scriptsize
\centering
\caption{comparison between SD POD-TPWL and FD method for case 1}\label{tab4}
\begin{spacing}{1.25}
\begin{tabular}{|c|c|c|c|c|c|}
\hline
- & Iterations & FOM & $J(\boldsymbol{\xi})$$\times$10$^{4}$ & $e_{obs}$ & $e_{\boldsymbol{\beta}}$\\
\hline
Initial model & - & - & 1.69 & 28.38  & 2.28 \\
SD POD-TPWL & 103 & 55 & 0.0160  & 3.35   & 0.68   \\
FD & 52 & 988 & 0.0153  & 3.28  & 0.72 \\
'True' model & - & - & 0.0068  & 2.12  & 0 \\
\hline
\end{tabular}
\end{spacing}
\end{table}

\begin{figure}[!h]
\centering\includegraphics[width=1.1\linewidth]{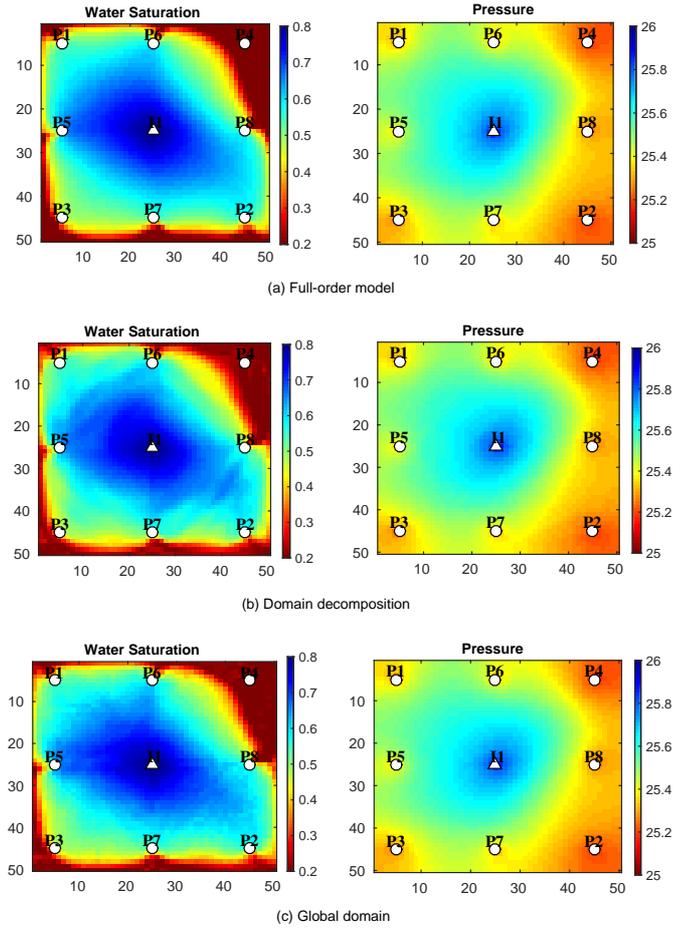}\\
\caption{Water saturation and pressure fields from the full-order model and from SD POD-TPWL and GD POD-TPWL based models for case 1.}\label{fig12}
\end{figure}

One of the key issues for the subdomain POD-TPWL is the implementation of the domain decomposition technique.
Our proposed subdomain POD-TPWL (SD POD-TPWL) can be easily generalized to the global domain POD-TPWL (GD POD-TPWL).
The differences between SD POD-TPWL and GD POD-TPWL are: a) model order reduction in global domain versus redcution in each subdomain separately; b) derivative estimation using RBF interpolation
in the global domain versus interpolation for each subdomain. As shown in Table \ref{tab3}, the total dimension of the reduced-order linear model is 18+122+51=191 for domain decomposition and 18+72+36=126
for the global domain. Table 5 shows the total number of the reduced variables in each subdomain and in the global domain. While the total sum of the
reduced variables in each subdomain is larger than that of the global domain, the number of reduced variables in each individual subdomain is relatively small.
Furthermore, these local reduced variables have the surprisingly abilities to accurately capture the flow dynamics, as suggested by Fig.\ref{fig12}. Fig.\ref{fig12}
shows the distribution of pressure and saturation at the final time. In this case, the reconstructions of the saturation and pressure field using a small
number of patterns in each subdomain are comparable with those of the global domain. In addition, as shown in Table.\ref{tab6}, both GD POD-TPWL and SD POD-TPWL can converge to satisfactory results.
The SD POD-TPWL needs 55 FOM simulations, while the GD POD-TPWL algorithm requires 73 FOM simulations (15 FOM simulations are run to collect the snapshots, 55 FOM simulations are used to construct
the initial reduced-order linear model in the first outer-loop, and the remaining 3 FOM simulations are used to reconstruct the reduced-order linear models in the following three outer-loops).
Therefore, compared to the global RBF interpolation, the proposed local RBF interpolation technique requires only a small number of reduced variables per subdomain and is much more computationally efficient.
If the dimension of the underlying model would be much larger, the GD POD-TPWL would result in a reduced-order linear model with a higher dimension and therefore more interpolation points would be required in the RBF scheme.
In the SD POD-TPWL algorithm this problem is avoided since for large-scale problems the dimension of the reduced-order linear model for the subdomain does not increase significantly; we only need to activate more subdomains.

\begin{table}[!h]
\footnotesize
\centering
\caption{The number of interpolation variables in each subdomain and global domain for case 1. $\Omega^{d}$ is the \textit{d}-th subdomain}\label{tab5}
\begin{spacing}{1.25}
\begin{tabular}{|c|c|c|}
\hline
\multicolumn{2}{|c|}{Domain Decomposition} & Global Domain \\
\hline
1 &  75=21($\Omega^{1}$)+19($\Omega^{2}$)+17($\Omega^{4}$)+18  & \multirow{9}{*}{126 =72 +36+18}  \\
2 &  98=21($\Omega^{1}$)+19($\Omega^{2}$)+17($\Omega^{3}$)+23($\Omega^{5}$)+18 & \multirow{9}{*}{}  \\
3 &  74=19($\Omega^{2}$)+17($\Omega^{3}$)+20($\Omega^{6}$)+18  & \multirow{9}{*}{}  \\
4 &  97=21($\Omega^{1}$)+17($\Omega^{4}$)+23($\Omega^{5}$)+18($\Omega^{7}$)+18 & \multirow{9}{*}{}  \\
5 &  118=19($\Omega^{2}$)+17($\Omega^{4}$)+23($\Omega^{5}$)+20($\Omega^{6}$)+21($\Omega^{8}$)+18 & \multirow{9}{*}{}   \\
6 &   95=17($\Omega^{3}$)+23($\Omega^{5}$)+20($\Omega^{6}$)+17($\Omega^{9}$)+18 & \multirow{9}{*}{}  \\
7 &   74=17($\Omega^{4}$)+18($\Omega^{7}$)+21($\Omega^{8}$)+18 & \multirow{9}{*}{}  \\
8 &   97=23($\Omega^{5}$)+18($\Omega^{7}$)+21($\Omega^{8}$)+17($\Omega^{9}$)+18 & \multirow{9}{*}{}  \\
9 &   76=20($\Omega^{6}$)+21($\Omega^{8}$)+17($\Omega^{9}$)+18 & \multirow{9}{*}{}  \\
\hline
\end{tabular}
\end{spacing}
\end{table}

For Case 1, history matching results using GD POD-TPWL are slightly better than those from the subdomain POD-TPWL, especially for the high-permeable zone, e.g, the red area in Fig.\ref{fig9}.
The water-front of the waterflooding process propagates forward quickly (as the blue area in Fig.\ref{fig12}) and therefore there are strong dynamic interactions within this area.
Our chosen domain decomposition may artificially cut off this inherent dynamic interaction between the east-south corner and the west-north corner.
A flow-informed domain decomposition technique may therefore be required to identify the relevant dynamic interactions,
especially for strongly heterogeneous reservoir models such as those based on strongly contrasting facies distributions or channels.

\begin{table}[!h]
\scriptsize
\centering
\caption{comparison between SD POD-TPWL and GD POD-TPWL for case 1}\label{tab6}
\begin{spacing}{1.25}
\begin{tabular}{|c|c|c|c|c|c|}
\hline
-  & FOM & $J(\boldsymbol{\xi})$$\times$10$^{4}$ & $e_{obs}$ & $e_{\boldsymbol{\beta}}$\\
\hline
Initial model & - & 1.69 & 28.38  & 2.28 \\
SD POD-TPWL & 55 & 0.0160  & 3.35   & 0.68   \\
GD POD-TPWL & 73 & 0.0140  & 3.21   & 0.61   \\
FD & 988 & 0.0153  & 3.28  & 0.72 \\
'True' model & - & 0.0068  & 2.12  & 0 \\
\hline
\end{tabular}
\end{spacing}
\end{table}


Solutions in our previous numerical experiments do not enable us to quantify the uncertainty of the permeability field and the predictions.
In general random maximum likelihood (RML) procedure \cite{oliver2008inverse} enables the assessment of the uncertainty by generating multiple 'samples' from the posterior distribution.
Each of these samples is a history matched realization, which also honors the data. Traditional 4D-Var or gradient-based history matching method obtains only one specific solution.
Additional solutions are obtained by repeatedly implementing the minimization process, but has a very highly computational cost.
The RML procedure can be efficiently implemented using our proposed reduced-order history matching algorithm.
When different background parameters are chosen to construct the reduced-order linear models, several valid solutions are obtained based on acceptable data misfits.
In this case, we choose 20 different background parameter sets to repeatedly implement our proposed adjoint-based reservoir history matching process.

\begin{figure*}[!h]
\centering\includegraphics[width=1\linewidth]{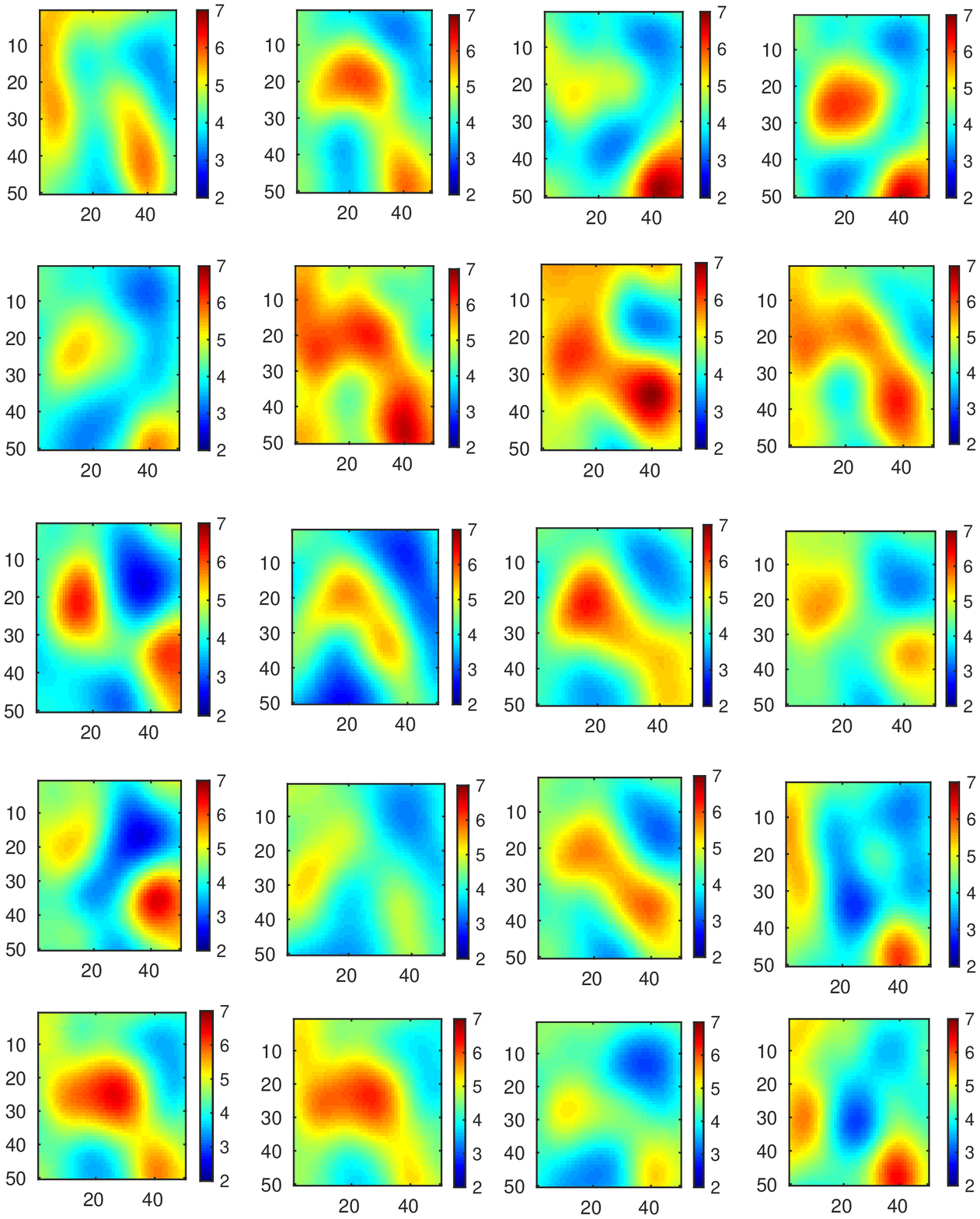}\\
\caption{Ensemble of 20 updated log-permeability fields for case 1.}\label{fig13}
\end{figure*}

\begin{figure}[!h]
\centering\includegraphics[width=1\linewidth]{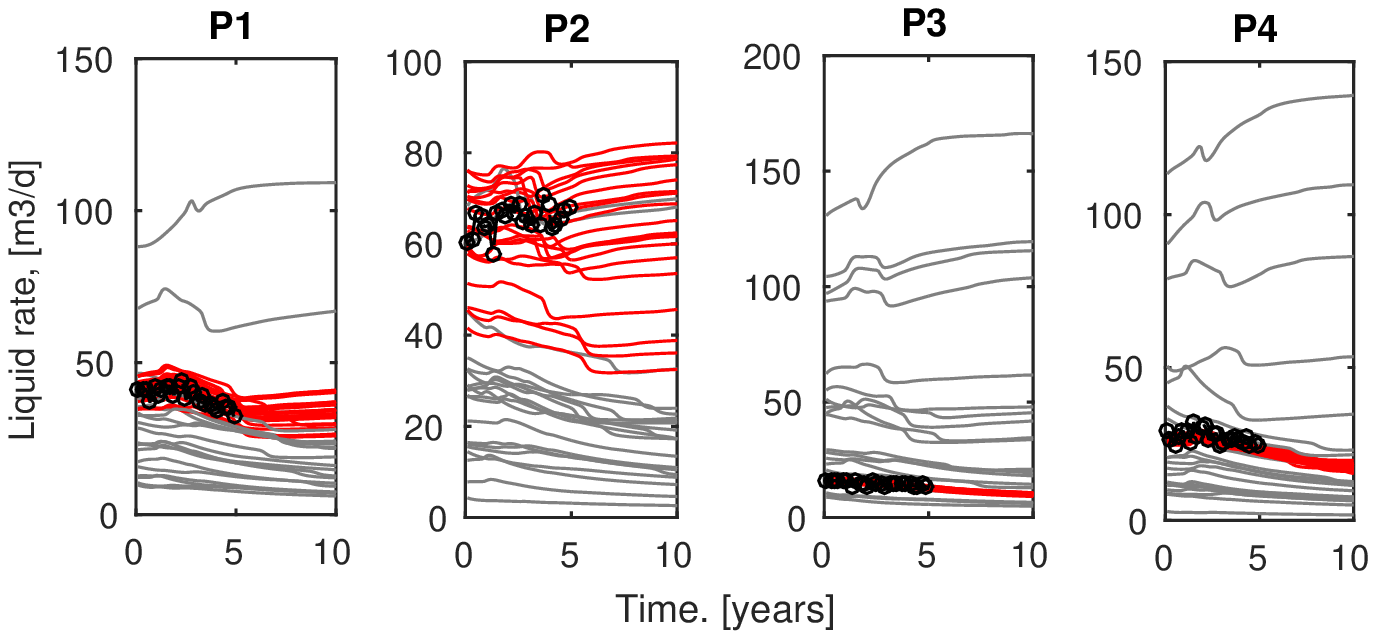}\\
\centering\includegraphics[width=1\linewidth]{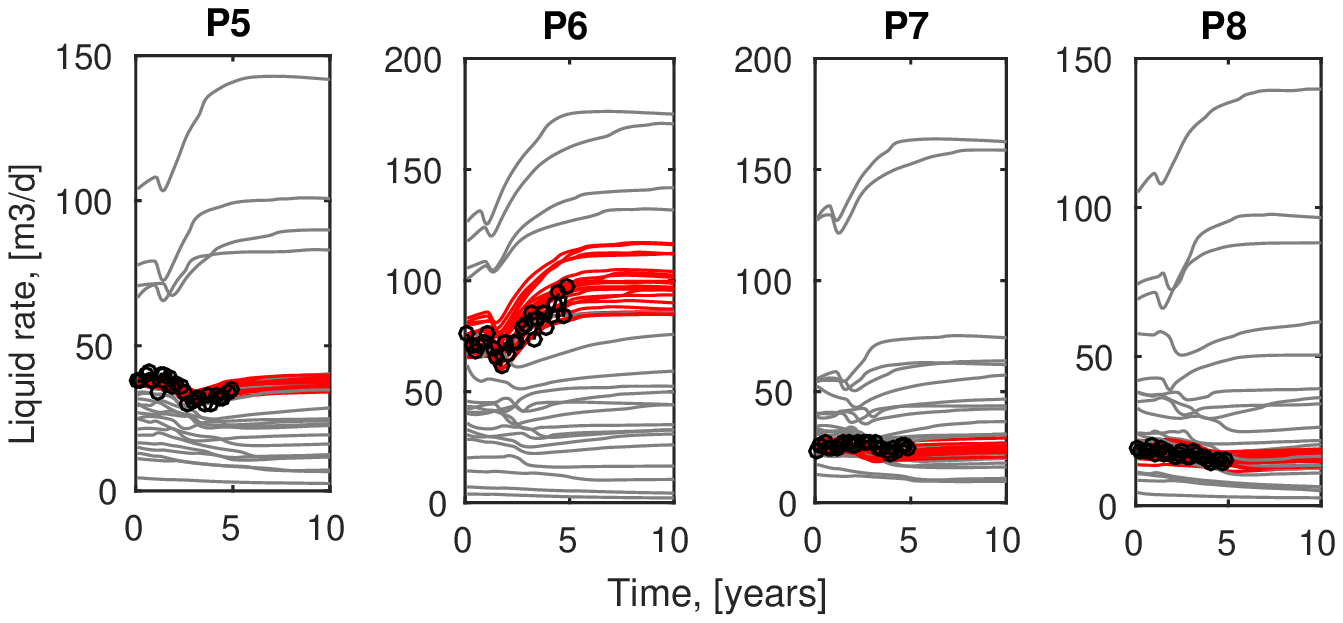}\\
\centering\includegraphics[width=1\linewidth]{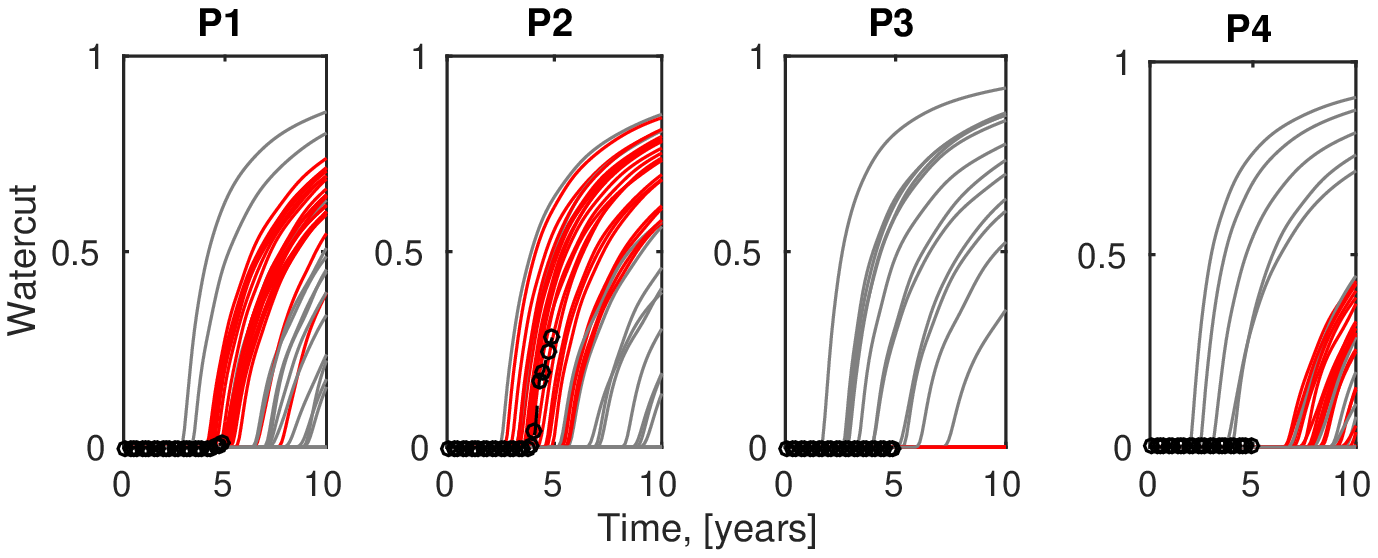}\\
\centering\includegraphics[width=1\linewidth]{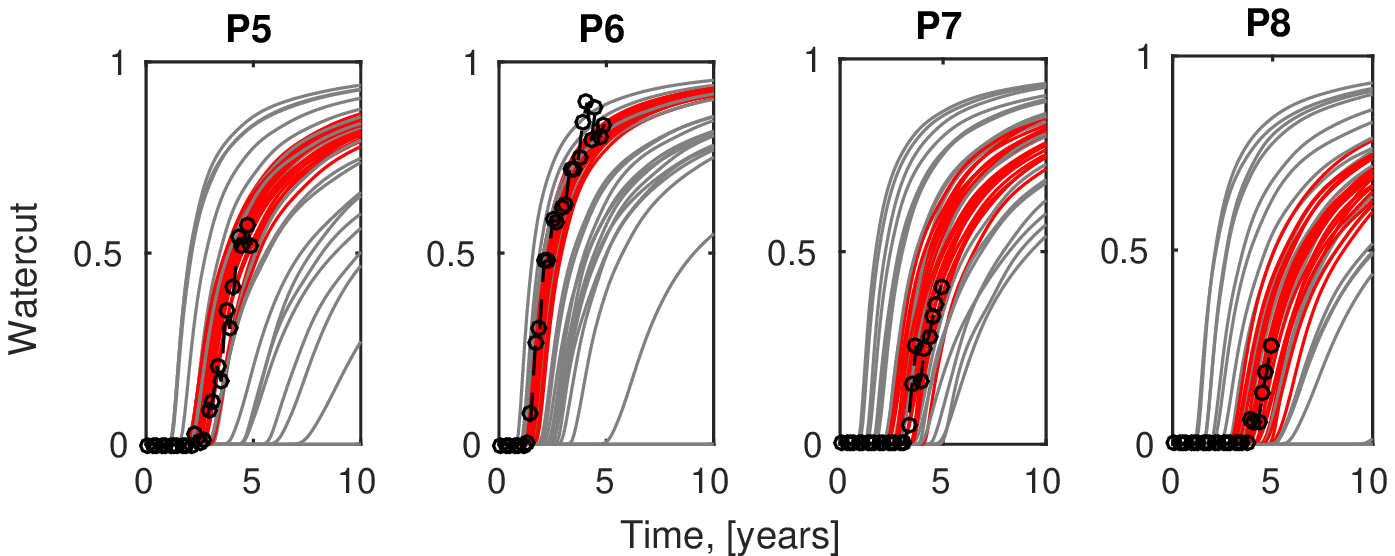} \\
\caption{ Ensemble of fluid rate and water-cut prediction for case 1. The gray lines represent the predictions from the 20 prior permeability realizations,
while the red lines represent the predictions from the corresponding 20 posterior permeability realizations calibrated using our method. The circles represent the noisy data.}\label{fig14}
\end{figure}

\begin{figure}[!h]
\centering\includegraphics[width=1\linewidth]{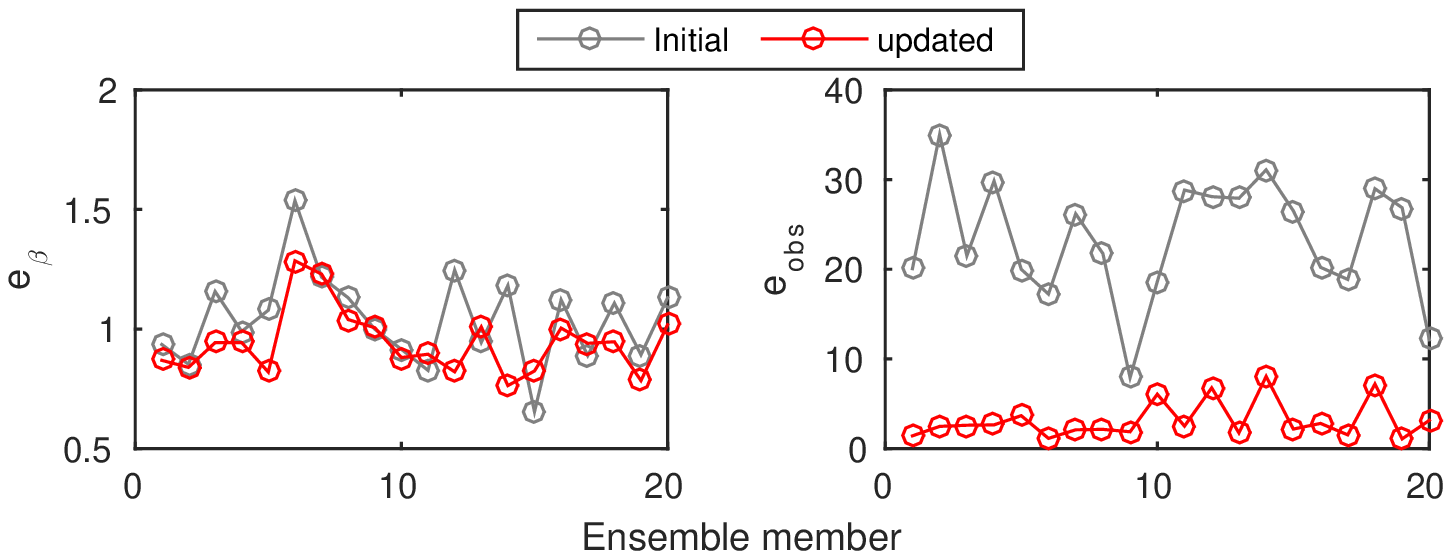}\\
\caption{Error analysis in terms of measurement and permeability misfits for the 20 solutions for case 1.}\label{fig15}
\end{figure}

Once the reduced model has been constructed the history matching can be efficiently repeated for different initial (background) models.
Fig.\ref{fig13} shows an ensemble of the posterior realizations (updated log-permeability field) using 20 different background parameters sets.
The main geological features, e.g, the high permeable area, are partly reconstructed in all of these 20 cases. Fig.\ref{fig14} and Fig.\ref{fig15} show an
ensemble of forecasts and the corresponding data misfits respectively using these 20 different initial and posterior models. Almost all of these 20 calibrated models produce
improved predictions of the fluid rate and WCT for all eight producers that are generally consistent with the data.
Thus, all of these 20 updated log-permeability fields can be regarded as acceptable solutions of the history matching problem.
The spread of the predictions from the posterior realizations is significantly decreased relative to the predictions from the prior realizations.
In the course of uncertainty quantification,  we randomly choose these 20 background parameter sets from the prerun 52 FOM simulations which is used to construct the reduced-order linear in the first outer-loop.
In addition, the analysis of the ensemble spreading enables us not to implement the additional outer-loops for updating the reduce-order linear model.
This means there is no need to run additional FOM simulations for the outer-loops, which makes our method very efficient.
Finally, in order to obtain these 20 solutions, the RML procedure totally requires 52 FOM simulations, including 15 FOM simulations for collecting snapshots,
37 FOM simulations for the initial construction of the reduced-order linear model and no additional outer-loop.

\subsection{Computational aspects}
This section discusses the computational aspects of our proposed adjoint-based reservoir history matching algorithm. The offline computational costs for subdomain POD-TPWL algorithm comprise (1) executing reparameterization using eigenvalue decomposition of the covariance matrix, (2) implementing model order reduction using POD in each subdomian, (3) conducting RBF interpolation and computing the derivative matrices. The cost of eigenvalue decomposition and POD is negligible for small models, while it will become significant for large-scale models. In our cases, the required number of FOM simulations is roughly 2-3 times the number of PCA coefficients, e.g, 54 simulations for the synthetic model, 113 (sampling within a small interval [-0.1, 0.1]) and 199 (sampling within a relative large interval [-1,1]) FOM simulations for the SAIGUP model, respectively. This process is code non-intrusive without the need of large programming effort. Besides, this process is also easily parallelized. Once available, the costs of running the reduced model are negligible. We should note that the gradient-based reservoir history matching  generally requires $O(10^{2}-10^{4})$ FOM simulations, thus, an offline cost of $O(10-10^{2})$ FOM simulations in these settings is attractive. For large-scale reservoir history matching, the main computational cost is dominated by the required number of FOM simulations. In our proposed method, most part of the FOM simulations is mainly in offline stage, which means that our method is easily implemented.

\section{Conclusions}
We have introduced a variational data assimilation method where the adjoint model of the original high-dimensional non-linear model is replaced by a subdomain reduced-order linear model.
Reparameterization and proper orthogonal decomposition techniques are used to simultaneously reduce the parameter space and reservoir model.
In order to avoid the need for simulator code access and modification and numerous full-order model simulations,
we integrated domain decomposition and radial basis function interpolation with trajectory piecewise linearization to form a new subdomain POD-TPWL algorithm.
The reduced-order linear model is easily incorporated into an adjoint-based parameter estimation procedure.
The use of domain decomposition allows for large-scale applications since the number of interpolation points required depends primarily on the number of the parameters and not on the dimension of the underlying full-order model.

We used the subdomain model-reduced adjoint-based history matching approach to calibrate the unknown permeability fields of a 2D synthetic model with noisy synthetic measurements.
The permeability field is parameterized using a KL-expansion resulting in a small number of permeability patterns that are used to represent the original gridblock permeability.
The reservoir domains is divided into 9 subdomains. In the numerical experiment, our methodology accurately reconstructs the 'true' permeability field and
shows similar results as more classic finite-difference based history matching. Our method also significantly improves the prediction of fluid rate and water breakthrough time of production wells.
Without any additional full-order model simulations, our approach efficiently generates an ensemble of models that all approximately match the observations.
For the cases studied in this paper, the number of full-order model simulations required for history matching is roughly 2-3 times the number of the number of global parameter patterns.

There are a number of aspects of the proposed methodology that could possibly be improved. It was observed that sampling strategy has to be chosen with care to obtain an efficient implementation.
Some diagnostics could possibly be devised to determine if and how many additional sampling points need to be generated.
We have chosen somewhat arbitrary decompositions of the global domain into subdomains. It may be beneficial to choose the subdomains based on information about the main dynamical patterns.
Since in reservoir applications these patterns are strongly affected by the placement of producers and injectors the subdomain decomposition could possibly be informed by the well lay-out.
In this paper we considered a global parameterization of the log-permeability field where the PCA patterns are defined over the entire domain.
From a computational point of view, a local parameterization where the parameters are defined in each subdomain separately is very attractive.
Since in this case parameters can be perturbed independent of each other and the effects of all these perturbations can be computed with very few full-order model simulations.
The local parameterization technique is the focus of our ongoing research. Also more complex history matching problem should be tested to shwo whether very promising results could be obtained.

\section*{Acknowledgment}
We thank the research funds by China Scholarship Council (CSC) and Delft University of Technology.












\end{document}